\documentclass[11pt,english]{smfart}
\usepackage[T1]{fontenc}

\usepackage[ethiop,english]{babel}
\usepackage{ethiop}
\usepackage{amssymb,url,xspace,smfthm}
\usepackage{MnSymbol}

\usepackage{hyperref}

 \newcommand{\ethi}{\selectlanguage{ethiop}}
\DeclareMathAlphabet\mathbb{U}{msb}{m}{n}

\usepackage{graphicx}
\usepackage[matrix, arrow,all,cmtip,color]{xy}

\newcommand{\bi}{\begin{itemize}}
\newcommand{\ei}{\end{itemize}}
\newcommand{\bd}{\begin{description}}
\newcommand{\ed}{\end{description}}

\newcommand{\bee}{\begin{enumerate}}
\newcommand{\eee}{\end{enumerate}}

\theoremstyle{plain}

\newtheorem{lemma}{Lemma}

\theoremstyle{definition}

\newtheorem{remark}{Remark}

\def\bqqq{\begin{quote}}
\def\eqqq{\end{quote}}

\def\lra{\longrightarrow}
\def\ra{\rightarrow}
\def\llrra{\leftrightarrow}

\def\rtt{\,\rightthreetimes\,}

\makeatletter
\newcommand{\xleftrightarrow}[2][]{\ext@arrow 3359\leftrightarrowfill@{#1}{#2}}
\newcommand{\xdashrightarrow}[2][]{\ext@arrow 0359\rightarrowfill@@{#1}{#2}}
\newcommand{\xdashleftarrow}[2][]{\ext@arrow 3095\leftarrowfill@@{#1}{#2}}
\newcommand{\xdashleftrightarrow}[2][]{\ext@arrow 3359\leftrightarrowfill@@{#1}{#2}}
\def\rightarrowfill@@{\arrowfill@@\relax\relbar\rightarrow}
\def\leftarrowfill@@{\arrowfill@@\leftarrow\relbar\relax}
\def\leftrightarrowfill@@{\arrowfill@@\leftarrow\relbar\rightarrow}
\def\arrowfill@@#1#2#3#4{%
  $\m@th\thickmuskip0mu\medmuskip\thickmuskip\thinmuskip\thickmuskip
   \relax#4#1
   \xleaders\hbox{$#4#2$}\hfill
   #3$%
}
\makeatother
\def\xra{\xrightarrow}
\def\xdra{\xdashrightarrow}

\def\ZZ{\Bbb Z}

\def\FFF{\mathfrak F}

\def\UUU{\mathfrak U}

\def\xra{\xrightarrow}

\def\rrt#1#2#3#4#5#6{\xymatrix{ {#1} \ar[r]^{} \ar@{->}[d]_{#2} & {#4} \ar[d]^{#5} \\ {#3}  \ar[r] \ar@{-->}[ur]^{}& {#6} }}

\def\NN{\Bbb N}
\def\RR{\Bbb R}

\def\XxK{X\times K}
\def\UxV{U\times V}
\def\EA{\exists\forall}
\def\AE{\forall\exists}

\def\FF{\mathcal F}

\def\NN{\mathbb N}

\def\oo{\infty} 
\def\id{{\text{id}}}
\def\dist{\text{dist}}

\def\Ord#1{\mathrm{Ord}_{<#1}}

\def\Filt{{\ethi\ethmath{wA}}\!\mathit{ilt}}
\def\sFilt{s{{\ethi\ethmath{wA}}}\!\mathit{ilt}}

\def\FFilt{{\ethi\ethmath{wE}}\!\text{\it ilt}}
\def\sFFilt{s{\ethi\ethmath{wE}}\!\text{\it ilt}}

\def\Fmilt{{\ethi\ethmath{qa}}\!\text{\it ilt}}
\def\sFmilt{s{\ethi\ethmath{qa}}\!\text{\it ilt}}

\def\FFilt{{\ethi\ethmath{wE}}\!\text{\it ilt}}
\def\sFFilt{s{\ethi\ethmath{wE}}\!\text{\it ilt}}

\def\Fmilt{{\ethi\ethmath{qa}}\!\text{\it ilt}}
\def\sFmilt{s{\ethi\ethmath{qa}}\!\text{\it ilt}}

\def\Ee{\iembE}

\def\iemb{{\ethi{\ethmath{nI}}}}

\def\iembE{{\ethi{\ethmath{hI}}}}

\def\ttt{{\ethi{\ethmath{pa}}}}

\def\mU{{\ethi{\ethmath{mi}}}}

\begin{document}
\selectlanguage{english}\catcode`\_=8
\date{}
\title[A naive  approach to tame topology]{
A naive diagram-chasing approach to formalisation of tame topology
}
\author[]{
notes by misha gavrilovich and konstantin pimenov\thanks
{A draft of a research proposal. 
Comments welcome at   
 $\tt{mi\!\!\!ishap\!\!\!p@sd\!\!\!df.org}$. 
I thank Martin Bays, Sergei Ivanov and Vladimir Sosnilo for discussions. {\tiny\href{http://mishap.sdf.org/mintsGE/}{http://mishap.sdf.org/mints/}} 
}
\\
{\sf\small in memoriam: evgenii shurygin}\\
{\small---------------------------------------------------------------------------------------------------------} \\
\tiny 
{\tiny ...~due to internal constraints on
    possible architectures of unknown to us~functional~"mental~structures".
 }\\ {\tiny
  Misha Gromov. Structures, Learning and Ergosystems: Chapters 1-4, 6.
 } 
}
\maketitle
\begin{abstract}
We rewrite classical topological definitions using the category-theoretic 
notation of arrows
and are thereby led to their concise reformulations in terms of
simplicial categories and 
orthogonality of morphisms, which we hope might be of use
in the formalisation of topology and in developing the tame topology of Grothendieck.  

Namely, we observe that topological and uniform spaces are simplicial objects 
in the same category, a 
category of filters, and 
that a number of elementary properties can be  
obtained by repeatedly passing
  to the left or right orthogonal (in the sense of Quillen model categories) 
  starting from a simple class of morphisms, often
  a single typical (counter)example appearing implicitly in the definition. 

Examples include the notions of: compact, discrete, connected, and
 totally disconnected spaces, dense image, induced topology, and separation axioms,
and, outside of topology, finite groups being
nilpotent, solvable, torsion-free, $p$-groups, and prime-to-$p$ groups;
 injective and projective
 modules; injective and surjective (homo)morphisms.

\end{abstract}
\setcounter{tocdepth}{4}
\tableofcontents

\section{
Introduction.
}
\subsection{Main ideas.}
In this note we rewrite several classical definitions and constructions in topology  in terms
of category theory and diagram chasing. We do so by
by first ``transcribing'' excerpts of [Bourbaki, General Topology] and [Engelking, Topology] 
by means of notation extensively using arrows, and then recognizing familiar patterns of 
standard category-theoretic constructions and diagram chasing arguments.

Arguably,
{\em we transcribe the ideas of Bourbaki into a language of category theory
appropriate to these ideas}, and our analysis of the text of Bourbaki
shows these ideas (but not notation) are implicit in Bourbaki and reflect their
logic (or perhaps their ergologic in the sense of 
\href{http://www.ihes.fr/~gromov/PDF/ergobrain.pdf}{[Gromov. Ergobrain; }
\href{http://www.ihes.fr/~gromov/PDF/ergo-cut-copyOct29.pdf}{Memorandum Ergo]}).

Doing so, we observe that a number of elementary textbook properties 
are obtained by taking the orthogonal (in the sense of Quillen lifting property) 
 to the simplest morphism-counterexample, and this 
leads to a concise syntax expressing these properties in two or three bytes
in which e.g. denseness, separation property Kolmogoroff/$T_0$,  compactness is expressed as
$$
  \begin{array}{ccccccc}
  \text{(dense image)}
     & & 
  \text{(Kolmogoroff/$T_0$)}
    & & 
  \text{(compact)}
  \\
  (\{ c \}\lra \{ o \ra c \})^l &\ \ \ \ &
 (\{x\leftrightarrow y \}\lra \{ x =y \})^r & \ \ \ \ \ &
 ((\{ \{o\}\lra \{o\ra c\}\}^{r})_{<5})^{lr}   
\end{array}
$$
this shows their \href{https://arxiv.org/abs/1301.0081}{Kolmogoroff complexity} is very low.

We also observe that the categories of topological spaces,  uniform spaces, and  simplicial sets 
are all, 
in a natural way, 
full subcategories of the same larger category, namely the simplicial category of filters;
coarse spaces of large scale metric geometry are also simplicial objects of a category of filters with different morphisms. 
This is, moreover, implicit in the definitions of a topological, uniform, and coarse space.

The exposition is in form of a story where we pretend to ``read off''
category-theoretic constructions from the text of excerpts of [Bourbaki] and
[Engelking] in a straightforward, unsophisticated, almost mechanical manner.
We hope word ``mechanical'' can be taken literally: we pretend to search for
correlations between the structure of allowed category-theoretic
diagram-chasing constructions and the text of arguments in topology, and hope
this search can be done by  a short program.  

No attempt is made to develop a theory or prove a theorem: our goal is 
to explain the process of transcribing by working out a few
examples in detail. In fact,
we think that understanding and formalising this process 
is a very interesting question.

This note is a research proposal suitable for a polymaths project: 
transcribing topological arguments into category theory involves rather independent tasks:
finding topological arguments worth transcribing and working out the precise meaning of category theoretic reformulations
are best suited for general topologists; spotting category theoretic patterns is best suited for category theorists; 
working out  formal syntax is best suited for logicians.

We hope our way of translating might of use in the formalisation of topology 
and suggests an approach to
the tame topology of Grothendieck.

\subsection{Contents.} In \S\ref{intro:sur}, as a warm-up and an example of our translation, we discuss the definition of surjection;
in \S\ref{yoga:ort}, we suggest the intuition that  orthogonality is category-theoretic {\em negation}. 
Appendix~\S\ref{surinj} gives a verbose exposition of the same ideas aimed at a student. 

In \S2.1 we start with a detailed translation of the definitions by Bourbaki of a dense subspace 
and a separation axiom of being Kolmogoroff/$T_0$ and show these definitions implicitly 
describe the simplest counterexamples involving spaces consisting of one or two points, and 
in fact require orthogonality to these counterexamples. Appendix~\S\ref{app:rtt-examples} and \S\ref{app:rtt-top}
gives more examples of properties defined by iterated orthogonals. Examples include the notions of: compact, discrete, connected, and
 totally disconnected spaces, dense image, induced topology, and separation axioms.
Appendix~\S\ref{app:top-notation} introduces a formal syntax and semantics which expresses these properties in several bytes
in both human- and a computer- readable form.
Outside of topology, examples in \S\ref{app:rtt-top} include finite groups being
nilpotent, solvable, torsion-free, $p$-groups, and prime-to-$p$ groups;
 injective and projective
 modules; injective, surjective.

Compactness is discussed in \S\ref{comp:ult} we reformulate the Bourbaki's definition in terms of convergence of ultrafilters
as an iterated orthogonal of the simplest counterexample. With help of this, we show in \S\ref{comp:m2} that 
there is a factorisation system corresponding to Stone-\v Cech compactification, and thus 
it is somewhat analogous to Axiom M2 $(cw)(f)$- and $(c)(wf)$- decompositions 
required in Quillen model categories.

In \S3, we ``transcribe'' the informal considerations in [Bourbaki, Introduction]. We ``read off'' from there in \S\ref{def:topoic-top} and \S\ref{met:filt}
that topological and uniform spaces are 2-dimensional simplicial objects 
in the same category, the category of filters. 
The discussion in \S\ref{exp:limits} of the notion of a limit of a filter $\FFF$ on a topological space $X$ 
leads to a reformulation in terms of the map \catcode`\_=8
$$(\FFF,\FFF\times \FFF,...)\lra (X\times X_{\ttt}, X\times X\times X_\ttt, ...)$$ 
from the object of Cartesian powers of $\FFF$ 
to the shift (d\'ecalage) of the simplicial object $(X_\ttt,X\times X_\ttt,...)$ corresponding to $X$. 
We end the section with a discussion in \S\ref{cwf-decompositions} of path spaces and cylinder objects in the category of topological spaces;
this also leads to constructions reminiscent of the shift (d\'ecalage).
Note that the d\'ecalage of a simplicial set is a model for the path space object of a topological space,
somewhat smaller than the usual model we discuss.\footnote{See~\href{https://ncatlab.org/nlab/show/decalage}{[nlab:d\'ecalage]} for a detailed discussion.}
 
In \S4, we formulate a number of open questions. Unfortunately, interesting open questions are rather vague and concern
the expressive power and formalisation of the new category theoretic, diagram-chasing way to talk about topology; to what extent the new language helps to avoid irrelevant set-theoretic details and counterexamples.
An important precise open question in this spirit is to define a model structure on the simplicial category of filters
compatible with a model structure on the full subcategory of topological spaces.

\subsection{Speculations.}
Does your brain (or your kitten's) have the lifting property (orthogonality), simplicial objects or diagram chasing built-in?
\S2 suggests a broader and more flexible context making
contemplating an experiment possible. Namely, some standard
arguments in point-set topology are computations with
category-theoretic (not always) commutative diagrams of finite categories (which happen to be preorders, 
or, equivalently, finite topological spaces) in the same way
that lifting properties define injection and
surjection. In that approach, the lifting property is viewed as a rule to add a new arrow,
a computational recipe to modify diagrams. 

Can one find an experiment
to check whether humans {\em subconsciously} use diagram chasing to reason about topology?

Does it appear implicitly in old original papers and books on point-set topology?

Is diagram chasing with preorders too complex to have evolved? Perhaps; but note the self-similarity:
preorders are categories as well, with the property that there is at most
one arrow between any two objects; in fact sometimes these categories are
thought of as $0$-categories. So essentially your computations are in
the category of (finite $0$-) categories.

Is it universal enough? Diagram chasing and point-set topology, arguably
a formalisation of ``nearness'',
is used as a matter of course in many arguments in mathematics.

Finally, isn't it all a bit too obvious? 
Curiously, in my experience it's a party topic people often get stuck on.
If asked, few if any can define a surjective or an injective map without words,
by a diagram,
or as a lifting property, even if given the opening sentence of \S\ref{app:sur-and-in} 
as a hint.
No textbooks seem to bother to mention these reformulations (why?).
An early version of \href{http://mishap.sdf.org/mints/Exercises_de_style_A_homotopy_theory_for_set_theory-I-II-IJM.pdf}{[Gavrilovich, Hasson]} states (*)${}_{\rtt}$ and (**)${}_{\rtt}$ of \S\ref{intro:sur} and \S\ref{app:sur-and-in}
as the simplest examples of lifting properties we were able to think up;
these examples were removed while preparing for publication.

\subsection{Surjection: an example}\label{intro:sur}  Let us now explain what we mean by translation. 
A map $f:X\lra Y$ is {\em surjective} iff it is left-orthogonal to the simplest non-surjective map $\emptyset\lra\{\bullet\}$, i.e.
$$\!\!\mathrm{(*)_{\rtt}}\ \ \ \ \ \ \ \ \ \ \ \ \ \ \ \ \ \  \emptyset\lra\{\bullet\} \rtt X \xra f Y$$
Recall that for morphisms $f:A\longrightarrow B$, $g:X\longrightarrow Y$ in a category, {\em a morphism $f$ has the left lifting property 
wrt a morphism $g$}, {\em  $f$ is (left) orthogonal to $g$}, and we write  $f
\,\rightthreetimes\,  g$ or $A \xra f B \,\rightthreetimes\, X \xra  g Y $, 
 iff for each $i:A\longrightarrow X$, $j:B\longrightarrow Y$ such that $ig=fj$ (``the square commutes''), 
there is $j':B\longrightarrow X$ such that $fj'=i$ and $j'g=j$ (``there is a diagonal 
making the diagram commute'').

With this definition, $\mathrm{(*)_{\rtt}}$ reads as
\begin{itemize}
\item[$\mathrm{(*)_{\text{words}}}$] 
{\em for each map $\{\bullet\}\xra y Y$, i.e. a point $y\in Y$, 
there is a map $\{\bullet\}\xra x X$, i.e. a point $x\in X$, such that $f\circ x=y$, i.e. $f(x)=y$}. 
\end{itemize}
This is the {\em text} of the usual definition of surjectivity of a function found in an elementary textbook. 
Conversely, we can read off $\mathrm{(*)_{\rtt}}$ from the text of the
definition of surjectively, by drawing the commutative diagram as we read
$\mathrm{(*)_{\text{words}}}$. 

It is this kind of direct, almost syntactic, relationship between 
the usual text and its category theoretic reformulation 
we are looking for in this paper. 
This is what we mean by saying {\em the reformulation  $\mathrm{(*)_{\rtt}}$
is implicit in the text  $\mathrm{(*)_{\text{words}}}$}.

For a property (class) \ensuremath{C} of arrows (morphisms) in a category,
define its {\em left} and {\em right orthogonals}, which we also call its 
{\em left} and {\em right negation} 
$$ C^l := \{ \ensuremath{f} :\text{ for each }g \in C\ \ensuremath{f} \,\rightthreetimes\,  \ensuremath{g} \} $$ 
$$ C^r := \{ \ensuremath{g} :\text{ for each }f \in C\ \ensuremath{f} \,\rightthreetimes\,  \ensuremath{g} \} $$ 
$$ C^{lr}:=(C^l)^r,\, C^{ll}:=(C^l)^l, ... $$

  Take $C=\{ \emptyset\longrightarrow \{*\} \}$ in $Top$. A calculation 
shows 
that a few of its iterated negations are meaningful: in $Top$,  
$C^r$ is the class of surjections (as we saw earlier), $C^{rr}$ is the class of subsets, 
$C^{rl}$ is the class of maps of form $A\longrightarrow A\sqcup D$, $D$ is discrete;
$\{\bullet\}\lra A$ is in $C^{rll}$ iff $A$ is connected;
$Y$ is totally disconnected iff $\{\bullet\}\xra y Y$ is in $C^{rllr}$ for each map  $\{\bullet\}\xra y Y$ (or, 
in other words, each point $y\in Y$).
$C^{l}$ 
 is the class of maps $A\longrightarrow B$ such that either $A\neq\emptyset$ or $A=B=\emptyset$. 
$C^{ll}=C^{lr}$ is the class of isomorphisms; 
$C^{lll}=C^{llr}=C^{lrl}=C^{lrr}=..$ is the class of all maps. 

Thus we see that already in this simplest case, taking iterated orthogonals (negation) produces 
several notions from a textbook, 
namely surjective, subset, discrete, connected, non-empty, and totally disconnected.

\subsection{Intuition/Yoga of orthogonality.}\label{yoga:ort} We suggest the following intuition/yoga is helpful.\footnote
{We were unable to find literature which explicitly describes this intuition, and will be thankful for any references which either discuss this intuition or list potential (counter)examples. 
}
\begin{itemize}

\item taking iterated orthogonals (negation) is a cheap way to automatically ``generate'' interesting notions;
a number of standard textbook notions are obtained in this way. 
We saw that taking iterated negations of the simplest map of topological spaces, $\{\}\lra\{\bullet\}$, 
generates $5$ classes worthy of being defined in a first year course of topology (surjective, subset, discrete, connected, 
non-empty and totally disconnected).  

\item it helps to think of {\em orthogonality  as  a category-theoretic (substitute for) negation};
taking orthogonal is perhaps the simplest way to define a class of morphisms 
without a property in a manner useful for a diagram chasing calculation. 

\item
often a morphism-counterexample  can be ``read off'' from the text of the definition of an elementary textbook property, 
and the property can be concisely reformulated as the orthogonal of the class
consisting of that counterexample.
\end{itemize}

\subsection{Intuition/Yoga of transcription.} We suggest the following intuition/yoga is helpful.
\begin{itemize}

\item
 ``transcribing'' the usual text of mathematical definitions and arguments    
by means of notation extensively using arrows sometimes makes it possible to recognise familiar patterns of
standard category-theoretic constructions and diagram chasing arguments.

\item orthogonality of morphisms often appears in this way, and so do simplicial objects

\item from the text of the definition of a topological property sometimes it is possible to ``read off'' 
a definition of a topology or a filter or a continuous function; it is worthwhile to try to interpret 
``for each open subsets there exists ...'' as a requirement that some function is continuous 

\ei

\section{Examples of translation. Orthogonality as negation.} 

\subsection{Dense subspaces and Kolmogoroff $T_0$ spaces.}
We shall now transcribe the definitions of {\em dense} and {\em Kolmogoroff $T_0$} spaces.

\subsubsection{``$A$ is a dense subset of $B$.''} 
By definition [Bourbaki, I\S1.6, Def.12], 
\newline\noindent\includegraphics[width=\linewidth]{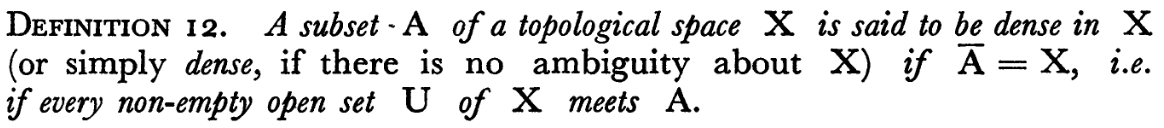}
 Let us transcribe this by means of the language of arrows.

{\sf $A$ is a subset of $B$}: this is an arrow $A\lra B$. 
(Note there is an alternative translation analogous to the used in the next item). 
{\sf Open subset}: An {\em open subset} of $B$
is an arrow $B\lra \{ U \ra U' \}$ ; here  $\{ U \ra U' \}$ denotes the topological space 
consisting of  one open point $U$ and one closed point $U'$; by the arrow $\ra$ we mean that
that $U'\in cl(U)$.
 {\sf Non-empty}: a subset $U$ is {\em empty} iff 
the arrow  $B\lra \{ U \ra U' \}$ factors as  $B\lra  \{ U' \}\lra \{ U \ra U' \}$ ;
here the map $ \{ U' \}\lra \{ U \ra U' \}$ is the obvious map sending $U'$ to $U'$. 
{\sf $U$ meets $A$}: $U\cap A =\emptyset$ iff the arrow $A\lra B\lra \{ U \ra U' \}$ factors as
$A\lra \{ U' \} \lra \{ U \ra U' \}$. 

Collecting above (Figure 1a), we see that 
a map $A\xra f B $ has dense image iff 
$$ A \xra f B \rtt  \{ U' \}\lra \{ U \ra U' \}$$ 

Note a little miracle: 
$\{ U' \}\lra \{ U \ra U' \}$ is the simplest map whose image isn't dense.
We'll see it happen again. 

\subsubsection{Kolmogoroff spaces, axiom $T_0$.} By definition [Bourbaki,I\S1, Ex.2b; p.117/122], 
\newline\noindent\includegraphics[width=\linewidth]{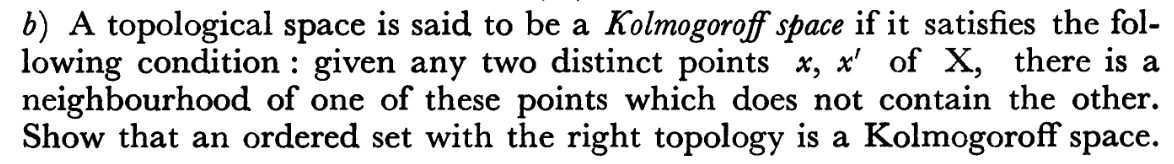}
Let us transcribe this. {\sf given any two ... points  $x$, $x'$ of $X$}: given a map $\{x,x'\}\xra f X$. 
{\sf two {\em distinct} points}:  the map $\{x,x'\}\xra f X$ does not factor through a single point, 
i.e.  $\{x,x'\}\lra  X$ does not factor as $\{x,x'\}\lra \{ x=x' \} \lra  X$. 
The negation of this sentence defines a topology on the set $\{x,x'\}$: indeed, 
 the antidiscrete topology on the set $\{x,x'\}$ is defined by the property that 
{\sf there is [no]  neighbourhood of one of these points which does not contain the other}
and we denote this space as  $\{x \llrra x'\}$. 
Now we note that the text implicitly defines 
the space $\{x \llrra x'\}$, and the only way to use it is to consider 
a map $\{x \llrra x'\}\xra f X$  instead of the map $\{x,x'\}\xra f X$.


Collecting above, we see that {\em a topological space $X$ is said to be a {\em Kolmogoroff} space 
iff any map $\{x \llrra x'\}\xra f X$ factors as $\{x \llrra x'\}\lra \{ x=x' \} \lra  X$.} 

Note another little miracle: it also reduces to  orthogonality of morphisms  
$$ \{x \llrra x'\}\lra \{ x=x' \} \rtt X\lra \{ x=x' \} $$
and $\{ x\llrra x' \}$ is the simplest non-Kolmogoroff space.

\subsubsection{Finite topological spaces as categories.}
Our notation  $\{ U' \}\lra \{ U \ra U' \}$ and $\{x \llrra x'\}\lra \{ x=x'
\}$ suggests that {\em we reformulated the two topological properties of being dense and Kolmogoroff 
in terms of diagram chasing in (finite)  categories}. And indeed, we may think
of finite topological spaces as categories and of continuous maps between them as {\em functors},
as follows; see Appendix~\ref{app:top-notation} for details and a definition of
our notation for finite topological spaces and maps between them.  

A {\em topological space} comes with a {\em specialisation preorder} on its points: for
points $x,y \in X$,  $x \leq y$ iff $y \in cl x$ ($y$ is in the {\em topological closure} of $x$).
 The resulting {\em preordered set} may be regarded as a {\em category} whose
{\em objects} are the points of ${X}$ and where there is a unique {\em morphism} $x{\ra}y$ iff $y \in cl x$.

For a {\em finite topological space} $X$, the specialisation preorder or
equivalently the corresponding category uniquely determines the space: a {\em
subset} of ${X}$ is {\em closed} iff it is
{\em downward closed}, or equivalently, 
it is a subcategory such that there are no morphisms going outside the subcategory.

The monotone maps (i.e. {\em functors}) are the {\em continuous maps} for this topology.

We denote a finite topological space by a list of the arrows (morphisms) in the
corresponding category; '$\leftrightarrow $' denotes an {\em isomorphism} and
'$=$' denotes the {\em identity morphism}.  An arrow between two such lists
denotes a {\em continuous map} (a functor) which sends each point to the
correspondingly labelled point, but possibly turning some morphisms into
identity morphisms, thus gluing some points.

\subsection{\label{comp:ult}Compactness via ultrafilters.}
We try to interpret the definition of compactness in [Bourbaki,I\S9.1, Def.1(C')]
in terms of arrows, or rather we try to rewrite it using the arrow notation
as much as possible. Doing so we shall see that   
this definition, in appropriate notation, condenses to
{\em a Hausdorff space $K$ is quasi-compact iff $K\lra\{\bullet\}$ is in 
$$  
 ((\{ \{o\}\lra \{o\ra c\}\}^{r})_{<5})^{lr},$$}
and that the latter expression almost appears in [Bourbaki, I\S10.2,Thm.1d] 
as a characterisation of the class of proper maps.

\subsubsection{Reading the definition of quasi-compactness.} We read the definition of quasi-compactness [Bourbaki,I\S9.1, Def.I];
we do not know how to read (C) and therefore we read reformulation (C$'$). 
\newline\noindent\includegraphics[width=\linewidth]{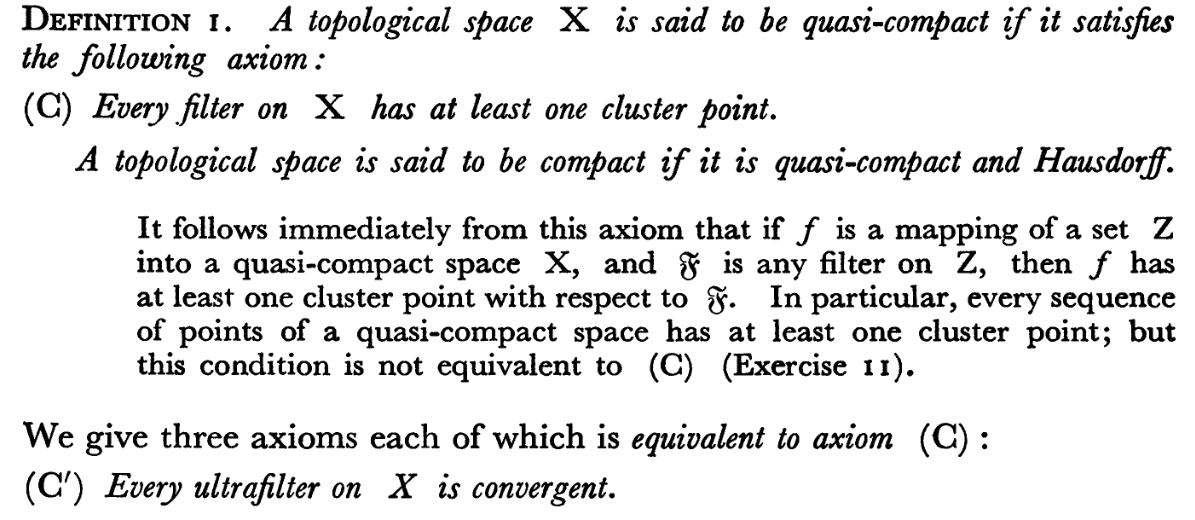}
 A space $K$ is {\em quasi-compact} iff each ultrafilter $\UUU$ on the set of points of 
$K$ converges,  
i.e. for each ultrafilter $\UUU$ on the set of points of
$K$ there is a point $x \in K$ such that each open neighbourhood of $x$ is $\UUU$-big. 
This contains a quantification over open subsets; this suggests to us that 
we should try to extract a definition of topology from the text and to interpret
the requirement as continuity of a certain map. 
{\sf each open neighbourhood of $x$ is $\UUU$-big} suggest we define a topology such that
an open subset is an $\UUU$-big open neighbourhood of some $x\in K$. This defines 
a topology on \def\oo{\text{``x''}} $K\sqcup \{\oo\}$: 
$$\{\,U \,:\, U\subset K\text{ is open}\}\cup \{\,U\cup\{\oo\} \,:\, U\subset K\text{ is open and }\UUU\text{-big}\}$$
 Denote the set equipped with this topology 
by $K\sqcup_{\UUU} \{\oo\}$. (Note [Bourbaki, I\S6.5, Definition 5, Example] define this space.)

Thus, in terms of arrows the definition becomes (see Figure 2a): 
$K$ is quasi-compact iff 
the identity map $K \xra {\id} K$ 
factors as $$K\lra K\sqcup_{\UUU} \{\oo\} \lra K$$ for each ultrafilter $\UUU$ on the set of points of $K$. 

Now note that Figure 2a is a particular case of orthogonality $K\lra K\sqcup_{\UUU} \{\oo\}  \rtt K\lra\{\bullet\}$, 
see Figure 2b where the map $K\lra K$ is arbitrary. Using orthogonals (negation), we express this by saying 
that  $K\lra K\sqcup_{\UUU} \{\oo\} \in \{  K\lra \{\bullet\}\}^l$. 
As usual,  we are tempted to define compactness as an orthogonal (negation) of a class (property) of morphisms, 
and therefore we check that  all maps of form $A\lra A\sqcup_{\UUU}\{\bullet\}$ lie
in this orthogonal $\{  K\lra \{\bullet\}\}^l$. Conversely, this also means that  $K\lra \{\bullet\}$, for $K$ quasi-compact, 
lies in the right orthogonal (negation)  $\{\,A\lra A\sqcup_{\UUU}\{\bullet\} \,:\, \UUU\text{ is an ultrafilter on a space }A\,\}^r$.

Summing up, we read Definition I as 
\begin{quote}
\noindent\includegraphics[width=\linewidth]{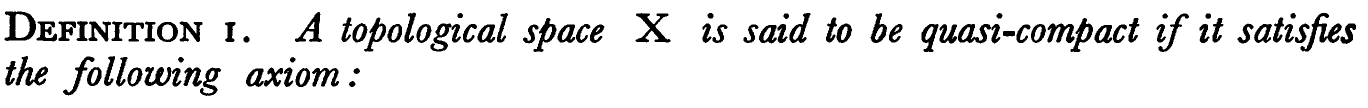}
$(\rm C')_{\rtt}$  $A\lra A\sqcup_{\UUU}\{\bullet\} \rtt X\lra \{\bullet\}$ for each ultrafilter $\UUU$ on each space $A$
\end{quote}

Note that there is another, more direct, way to read off the lifting property
from the remark in the proof of $(C)\implies (C')$:
\vskip3pt
\noindent\includegraphics[width=\linewidth]{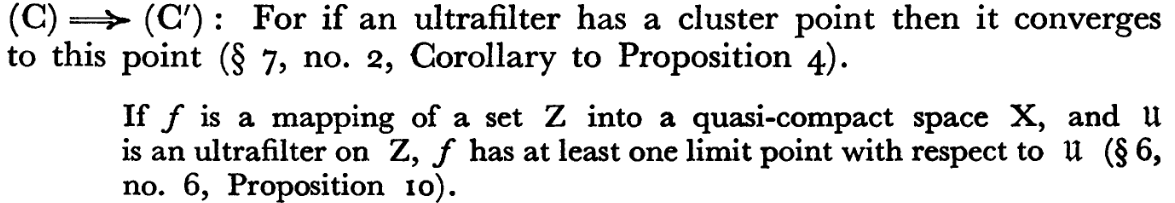}
%
%
In terms of arrows, this reformulation is {\em precisely} the lifting property
$$Z\lra Z\sqcup_{\UUU} \{\oo\}  \rtt X\lra\{\bullet\}$$
We'd like to view the fact that Bourbaki chooses to formulate explicitly {\em precisely} a lifting property 
 immediately following a key definition
as evidence that Bourbaki is implicitly doing category theoretic reasoning.

\subsubsection{Proper maps.} If we were to think that 
[Bourbaki, General Topology] does implicitly uses category theoretic reasoning and orthogonality, 
we'd hope to find there the definition of the class 
$$\{\,A\lra A\sqcup_{\UUU}\{\bullet\} \,:\, \UUU\text{ is an ultrafilter on a space }A\,\}^r$$
And indeed, this is how Bourbaki characterises the class of proper maps in 
[Bourbaki, General Topology, I\S10.2,Th.1(d)]
(cf.~Figure~2d), almost exactly.
We see this as evidence that Bourbaki does indeed use category theoretic reasoning, or perhaps
as an explanation of what do we mean by saying so.  

Note we might have started our translation with this characterisation of proper maps 
in terms of ultrafilters [Bourbaki, General Topology, I\S10.2,Th.1(d)], and 
we'd then arrive at Figure~2d directly.
\newline\noindent \includegraphics[width=\linewidth]{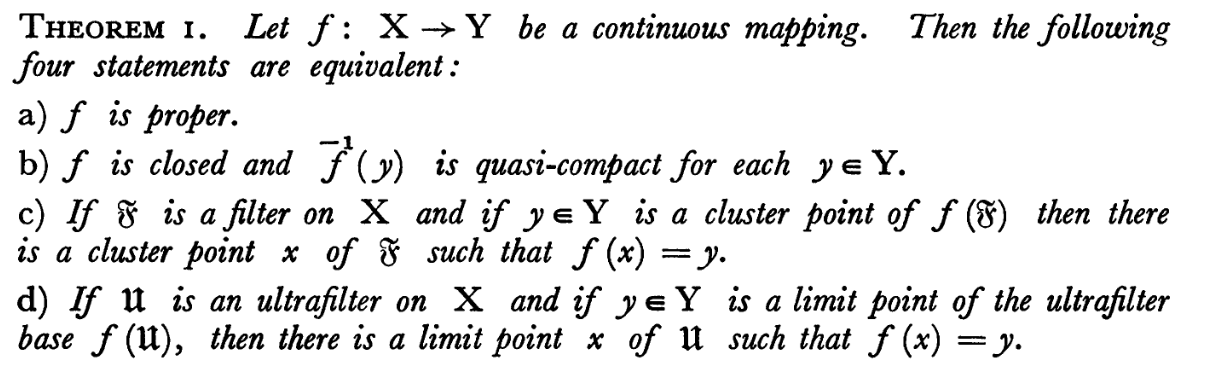}

However, this reformulation is unsatisfactory for us: 
it uses non-elementary, infinitary, set-theoretic
notion of ultrafilters which we do not know how to manipulate category-theoretically..

We'd like to have a definition which relies on maps between finite spaces.

An  argument similar to a linear algebra about dual vector spaces gives the following.
 For any class $C$ of maps we have that  $C^l=C^{lrl}$ and  $C^r=C^{rlr}$ 
and $C_1\subset C_2$ implies $C_1^l\supset C_2^l$ and $C_1^r \supset C_2^r$. This implies  
$P^{lr}\subset C^{rlr}=C^{r}$ whenever $P\subset C^r$.

Take $P$ to be some class of proper maps between finite spaces. By above we see that
$P^{lr}$ is a subclass of the class of proper maps. We want to take $P$ to be
large enough so that $P^{lr}$ is the whole class of  proper maps.
And indeed, we find that a classical theorem in general topology 
tells us we can do so, at least if we only care about spaces satisfying
separation axioms. Moreover, we will see it is enough to take $P$ to consist of 
the following 
maps between spaces of size at most 3:
$$\{\ 
{ \{\boldsymbol{B_1}\leftarrow O\rightarrow \boldsymbol{B_2}\} } \lra  {\{\bullet\} } \ ; \
 \{U\} \lra  { \{ U \ra U' \} } \ ; \
{ \{x\llrra y\} } \lra  \{x=y\} \ ; \
{ \{o \ra c \} } \lra  \{o=c\} \ \} 
$$

\subsubsection{Reducing to finite spaces.} 
Now we are back translating; we ignore the considerations of the previous subsubsection 
which give us a rather good idea of what we would get as the result of translation.

Reduction to finite spaces is provided by Smirnov-Vulikh-Taimanov theorem in the form by 
[Engelking, 3.2.1,p.136] (``compact'' below stands for ``compact Hausdorff''):\vskip3pt
\noindent \includegraphics[width=\linewidth]{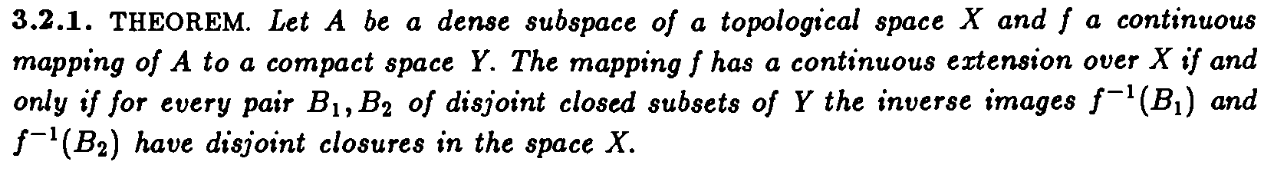}

Let us transcribe this. We are given $A\xra i X$ and  $A\xra f Y$. 
{\sf The mapping \ensuremath{f} has a continuous extension over \ensuremath{X}}
iff the arrow $A\xra f Y$ factors via $A\xra i X$ (cf.~Figure 2f). 
a {\sf   pair $B_1,B_2$ of disjoint closed subsets of \ensuremath{Y}} is an arrow 
$Y\lra  \{\boldsymbol{B_1}\leftarrow O \ra \boldsymbol{B_2}\}$
where $\{\boldsymbol{B_1}\leftarrow O \ra \boldsymbol{B_2}\}$ is the space with 
one open point denoted by $O$ and two closed points
denoted by $\boldsymbol{B_1}$ and $\boldsymbol{B_2}$.
{\sf  the inverse images $f^{-1}(B_1)$
and $f^{-1}(B_2)$ have disjoint closures in the space $X$} says 
the arrow $A\xra f Y \lra \{\boldsymbol{B_1}\leftarrow O \ra \boldsymbol{B_2}\}$ factors as $A \xra i X
\lra  \{\boldsymbol{B_1}\leftarrow O \ra \boldsymbol{B_2}\}$ (cf.~Figure 2g).

Now we need to define the class of  dense subspaces. We do so by 
taking left negations (orthogonals) of the simplest archetypal examples of 
maps with non-dense image, a non-injective map, and a map 
topology on the domain is not induced from the target. 


\begin{quote} 3.2.1. THEOREM. 
Let $Y$ be Hausdorff compact and let $A\xra i X$ satisfy (cf.~Figure 2(ijk))
\begin{enumerate}
\item[(i)]  (dense) $A\xra i X \rtt  \{U\} \lra  { \{ U \ra U' \} } $
\item[(ii)] (injective)  $A\xra i X \rtt { \{x\llrra y\} } \lra  \{x=y\}$
\item[(iii)] (induced topology) $A\xra i X \rtt { \{o \ra c \} } \lra  \{o=c\} $
\end{enumerate}
Then the properties of $A\xra f Y$ defined by Figure 2(f) and Figure 2(g) are equivalent.
\end{quote}

This implies that, for Hausdorff compact $Y$, 
items   3.2.1(i-iii) and $ A\xra i X \rtt \{\boldsymbol{B_1}\leftarrow O \ra \boldsymbol{B_2}\}\lra \{\boldsymbol{B_1} = O= \boldsymbol{B_2}\}$ imply that $ A\xra i X \rtt Y\lra \{\bullet\}$.

Further, note that if $X=A\sqcup\{\oo\}$ is obtained from $A$ by adjoining a single closed non-open point, then 
$$ A\xra i X \rtt \{\boldsymbol{B_1}\leftarrow O \ra \boldsymbol{B_2}\}\lra \{\boldsymbol{B_1} = O= \boldsymbol{B_2}\}$$
iff there exists an ultrafilter $\mathfrak U$ such that $A \xra i X$ is of form $A \lra A\sqcup_{\mathfrak U}\{\oo\}$.

This implies that maps of form  $A \lra A\sqcup_{\mathfrak U}\{\oo\}$ are in $P^l$ and, finally, 
that a Hausdorff space $K$ is quasi-compact iff $K\lra\{\bullet\}$ is in $P^{lr}$ where
$P$ is $$
\{\ 
{ \{\boldsymbol{B_1}\leftarrow O\rightarrow \boldsymbol{B_2}\} } \lra  {\{\bullet\} } \ ; \
 \{U\} \lra  { \{ U \ra U' \} } \ ; \
{ \{x\llrra y\} } \lra  \{x=y\} \ ; \
{ \{o \ra c \} } \lra  \{o=c\} \ \} 
$$

\subsubsection{The simplest counterexample negated three times.} 
Note that all maps between finite spaces mentioned in the preceeding subsubsection 
are closed, hence proper by  [Bourbaki, I\S10.2,Thm.1b].

A verification shows that, for $Y$ and $Z$ finite, the map  $Y\xra g Z $ is closed iff 
$$\{o\}\lra \{o \ra c \}  \rtt Y\xra g Z$$

Denote by 
$(\{ \{o\}\lra \{o\ra c\}\}^{r})_{<5}$  the subclass of $\{ \{o\}\lra \{o\ra c\}\}^{r}$ 
consisting of maps between spaces of size at most $4$. 

Considerations above could be summarized by:
\begin{itemize}
\item
a Hausdorff space $K$ is quasi-compact iff 
\begin{itemize}
\item[] $K\lra\{\bullet\}$ is in $  
 ((\{ \{o\}\lra \{o\ra c\}\}^{r})_{<5})^{lr}$.
\end{itemize}
\item every map in $((\{ \{o\}\lra \{o\ra c\}\}^{r})_{<5})^{lr}$ is proper
\end{itemize}

And we conjecture that the latter is in fact the class of all proper maps.

\begin{figure}
\small
\begin{center}

$ (a)\ \xymatrix{ K \ar[r]|{\text{id}} \ar@{->}[d] & K \ar[d] \\ { K\cup_{\FF}\{\text{``x''}\} } \ar[r] \ar@{-->}[ur]|{{\text{``x''}\mapsto x}}& \{\bullet\} }$
$\ \ \ \ (b)\ \xymatrix{ K \ar[r] \ar@{->}[d] & K \ar[d] \\ { K\cup_{\FF}\{\text{``x''}\} } \ar[r] \ar@{-->}[ur]  & \{\bullet\} }$
$\ \ \ \  (c)\ \xymatrix{ A \ar[r] \ar@{->}[d] & A \ar[d] \\ { A\cup_{\FF}\{\text{``x''}\} } \ar[r] \ar@{-->}[ur]  & \{\bullet\} }$

$ (d)\ \xymatrix{ X \ar[r]|{\text{id}} \ar@{->}[d] & X \ar[d]|{f} \\ { X\cup_{\mathfrak U}\{\text{``x''}\} } \ar[r] \ar@{-->}[ur]  & Y  }$
$\ \ \ \ \ (e)\ \xymatrix{ A \ar[r] \ar@{->}[d] & X \ar[d]|{g} \\ { A\cup_{\mathfrak U}\{\text{``x''}\} } \ar[r] \ar@{-->}[ur]  & Y  }$

$ (f)\ \xymatrix{ A \ar[r]|f \ar@{->}[d] & Y \ar[d]|{g} \\ { X } \ar[r] \ar@{-->}[ur]  & \{\bullet\}  }$
$\ \ \ \ \ (g)\ \xymatrix{ A \ar[r]|f \ar@{->}[d] & Y  \ar@{->}[r] & { \{\boldsymbol{B_1}\leftarrow O\rightarrow \boldsymbol{B_2}\} }  \\ { X }  \ar@{-->}[urr]  & { } & { } }$
$\ \ \ \  (h)\ \xymatrix{ A \ar[r] \ar@{->}[d] & { \{\boldsymbol{B_1}\leftarrow O\rightarrow \boldsymbol{B_2}\} } \ar[d] \\ { X } \ar[r] \ar@{-->}[ur]  & {\{\bullet\} }  }$

$ (i)\ \ \xymatrix{ A \ar[r] \ar@{->}[d] & { \{U\} } \ar[d] \\ { X } \ar[r] \ar@{-->}[ur]  & { \{ U \ra U' \} }  }$
$\ \ \ \ \ \ (j)\ \ \xymatrix{ A \ar[r] \ar@{->}[d] & { \{x\llrra y\} } \ar[d] \\ { X } \ar[r] \ar@{-->}[ur]  & \{x=y\}  }$
$\ \ \ \ \ \ (k)\ \ \xymatrix{ A \ar[r] \ar@{->}[d] & { \{o \ra c \} } \ar[d] \\ { X } \ar[r] \ar@{-->}[ur]  & \{o=c\}  }$

$ (l)\ \xymatrix{ {\{o\} } \ar[r] \ar@{->}[d] & X \ar[d] \\ { \{o \ra c \} } \ar[r] \ar@{-->}[ur]  & { Y }  }$

\end{center}
\caption{\label{fig1}\small
These are equivalent reformulations of quasi-compactness of spaces and its generalisation to maps, that of properness of maps. 
 (a) the identity map $K \xra {\id} K$
factors as $K\lra K\cup_{\FF} \{\oo\} \lra K$ 
 (b) this is also equivalent to $K$ being quasi-compact (we no longer require the arrow $K\lra K$ to be identity)
 (c) and in fact quasi-compact spaces are orthogonal to maps associated with ultrafilters 
 (d) $X\xra f  Y$ is proper, i.e. {\sf d) If $\mathfrak U$ is an ultrafilter on $X$ and if $y \in Y$ is a limit point of the ultrafilter
base $f (U)$, then there is a limit point $x$ of $\mathfrak U$ such that $f (x) = y$.}  [Bourbaki, General Topology, I\S10.2,Th.1(d)]
 (e) this is also equivalent to $X\xra f  Y$ is proper, i.e. this holds for each ultrafilter $\mathfrak U$ on each space $A$ 
 (f)  The mapping \ensuremath{f} has a continuous extension over \ensuremath{X}
 (h)  for every pair $B_1,B_2$ of disjoint closed subsets of \ensuremath{Y} the inverse images $f^{-1}(B_1)$
and $f^{-1}(B_2)$ have disjoint closures in the space $X$
 (i) the image of $A$ is dense in $B$
 (j) the map $A\lra B$ is injective
 (k) the topology on $A$ is induced from $B$
 (l) for $X$ and $Y$ finite, this means that the map $X\lra Y$ is closed, or, equivalently, proper
} \end{figure}

\subsubsection{\label{comp:m2}Compactification as factorisation system/M2-decomposition} 
By a simple diagram chasing argument,\footnote{\label{foot:bous}
See Thm.~3.1 of \href{https://core.ac.uk/download/pdf/82479252.pdf}{[Bousfield, Constructions of factorization systems in categories]}
for details of such an argument and assumptions which are enough to make it work. However, note that his definitions are somewhat different from ours: unlike us, he considers the {\em unique} lifting property, cf.\S2 [ibid.].} 
\footnote{ 
See \href{http://www.maths.uwc.ac.za/~dholgate/Papers/DBHPhd.zip}{[Holgate,PhD,2.1(Perfect Maps)]} and references therein
for examples of factorisation systems related to Stone-\' Cech decomposition and proper maps. 
Note [Holgate] says ``perfect'' instead of ``proper'', as is common in topology.
}
each morphism $X\lra Y$ decomposes as either
$X\xra{(P)^{rl}}\cdot\xra{(P)^r} Y $ and  $X\xra{(P)^{l}}\cdot\xra{(P)^{lr}} $ 
whenever $P$ is a class of morphisms and the underlying category has enough limits and colimits.

We shall now see that Stone-\v{C}ech compactification is an example of such a decomposition 
when $P$ is the class of proper maps
and is thus somewhat 
analogous to the  $(cw)(f)$- and $(c)(wf)$-decomposition 
required by Axiom M2 of Quillen model categories. 

 Almost this observation is mentioned explicitly in \href{https://core.ac.uk/download/pdf/82479252.pdf}{[Bousfield, Constructions of factorization systems in categories]}\footnote{
In our notation $\mathcal M (\mathbf E_1)$ is almost $(\mathbf E_1)^r$ but not quite: 
$\mathcal M (\mathbf E_1)$ is the right orthogonal ($\rtt$-negation) with respect to the {\em unique} lifting property; [7]  is [S. MacLane, Categories for the Working Mathematician (Springer-Verlag, New York, 1971)].
}: 
\newline\noindent\includegraphics[width=\linewidth]{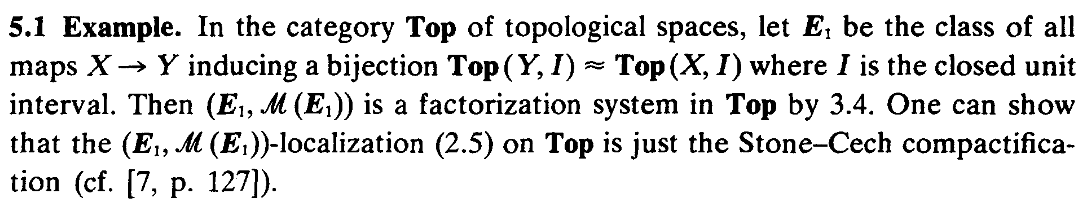}
\vskip-7pt
\noindent We shall find this observation by transcribing  [Engelking, Theorem 3.6.1, p.173] by means of diagram chasing. 
In fact, corollaries [Engelking, 3.6.2-3.6.9] could also be seen in a diagram chasing way; we only reformulate Corollary~3.6.3.
\begin{figure}
\noindent\includegraphics[width=\linewidth]{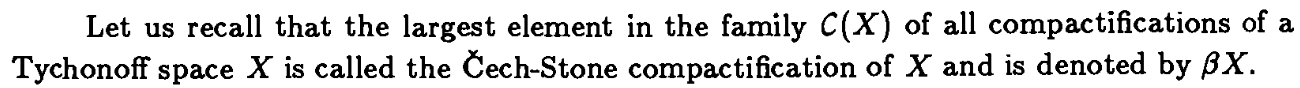}
\newline\noindent\includegraphics[width=\linewidth]{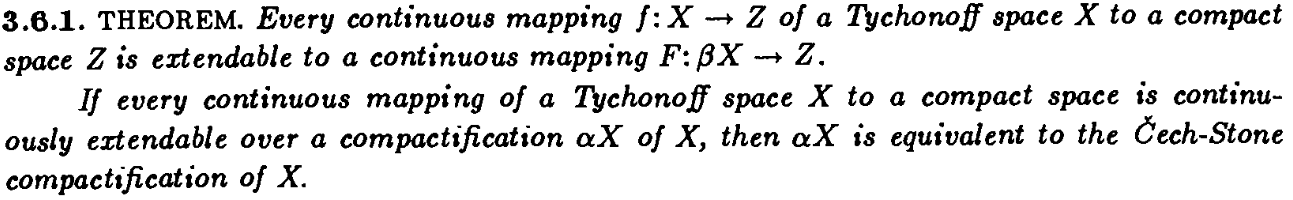}
\newline\noindent\includegraphics[width=\linewidth]{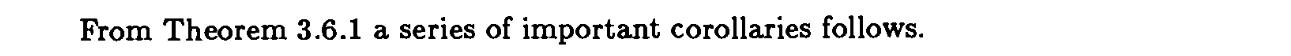}
\newline\noindent\includegraphics[width=\linewidth]{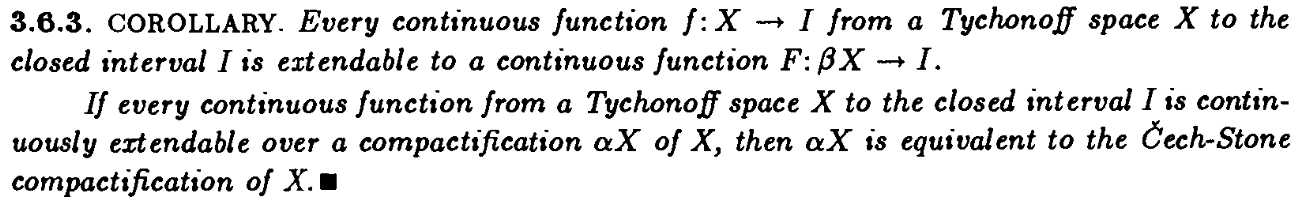}
\vskip5pt\small
\begin{center}
$ (a)\ \xymatrix{ X \ar[r]|{\text{f}} \ar@{->}[d] & Z \ar[d]|{(P)} \\ { \beta X } \ar[r]|{(P)} \ar@{-->}[ur]& \{\bullet\} }$
\ \ \ \ $ (b)\ \xymatrix{ X \ar[r]|{\forall} \ar@{->}[d] & {\beta X} \ar[d]|{(P)} \\ { \alpha X } \ar[r]|{(P)} \ar@{-->}[ur]& \{\bullet\} } $
\ \ \ \ $ (b')\ \xymatrix{ X \ar[r]|{(P)^l} \ar@{->}[d]|{(P)^l} & {\cdot} \ar[d]|{(P)^{lr}} \\ {\cdot} \ar[r]|{(P)^{lr}} \ar@{-->}[ur]|{(iso)}& \{\bullet\} } $

$ (c)\ \ \ X\xra{(P^l)} \beta X \xra {(P)} \{\bullet\}$ 
\ \ \ \ $ (d)\ \ \ \ X\xdra{(P^l)} Y \xdra {(P)^{lr}} Z$\\
$ (e) \ \ \ \ X\xra{(P^l)} \beta X \rtt [0,1]\lra\{\bullet\}$\\ 
 $(f)\ \ \ \ \ \  X\lra \alpha X \xra {(P)} \{\bullet\}$ and $ X\xra{(P^l)} \alpha X \rtt [0,1]\lra\{\bullet\}$ implies $\alpha X = \beta X$
\end{center}
\caption{\label{fig2}\small
A diagram chasing reformulation of [Engelking, Theorem 3.6.1, p.173]. 
(a) {\sf Every continuous mapping $f: X \lra Z$ of a Tychonoff space $X$ to a compact
space Z is extendable to a continuous mapping $F: \beta X \lra Z$.} 
(b) {\sf If every continuous mapping of a Tychonoff space X to a compact space is  
continuously extendable over a compactification $\alpha X$ of $X$, then $\alpha X$
 is equivalent to the Cech-Stone
compactification of $X$.} This is reformulated as follows: if diagram $(b)$ holds, 
then the diagonal map $\alpha X\lra \beta X$ can be chosen to be an isomorphism.
(b$'$) this is an analogue of (b) formulated in terms of category theory as uniqueness of $\cdot\xra{(P)^l}\cdot\xra{(P)^{lr}}\cdot$ decomposition; 
the diagonal arrow exists because $(P)^l\rtt (P)^{lr}$ 
and thus we require it to be an isomorphism. 
(c) Both diagrams above can summarized as: there exists a unique decomposition of this form. 
(d) Further, this is implied by an analogue of Axiom M2 $(cw)(f)$- and $(c)$$(wf)$-decomposition 
of model categories: each morphism $X\lra Z$ decomposes as 
$ X\xdra{(P^l)} Y \xdra {(P)^{lr}} Z$
(e) 
{\sf Every continuous function $f: X \lra X$ from a Tychonoff space $X$ to the
closed interval $I$ is extendable to a continuous function $F: \beta X\lra I$.}
(f) {\sf If every continuous function from a Tychonoff space $X$ to the closed interval $I$ is  
continuously extendable over a compactification $\alpha X$ of $X$, then $\alpha X$ is equivalent to the \v Cech-Stone
compactification of $X$.} Note the conclusion $\alpha X=\beta X$ is stated somewhat imprecisely; we rather need 
to say that morphisms $X\lra \alpha X$ and $X\lra \beta X$ are the same.
}\end{figure}
Let us transcribe this by means of the notation of arrows. We will deliberately ignore the separability assumptions that
$X$ is Tychonoff and $\beta X$ and $Z$ are assumed to be Hausdorff. 

Let $P$ be the class of proper maps. 
Figure 2ab represent the statement of the theorem. 
Figure 2a suggests that the compactification map $X\lra \beta X$ is in the class $(P)^l$;
Figure 2b suggests that there is a unique decomposition 
$X\xra{(P^l)} X' \xra {(P)^{lr}} \{\bullet\}$. 

And indeed, this is implied by a simple diagram chasing argument.
Uniqueness follows from orthogonality of $(P)^l$ and $(P)^r$. 
The decomposition is constructed by an argument which looks roughly as follows:\footnote{For details see footnote ${}^{(\ref{foot:bous})}$.}
consider all the decompositions of form $X\xra{(P^l)} X'\lra Y$ and 
take the pushout $X\xra{(P)^l} X_l$ of all the maps $X\xra{(P^l)} X'$ appearing in the decompositions of this form.
The map belongs to  $(P)^l$ because left orthogonals are closed under pushouts, 
By the universality property of pushouts you obtain a decomposition $X\xra{(P)^l} X_l\lra Y$
and a diagram chasing argument based on the definition of pushout and orthogonality properties of $(P)^l$ and $(P)^{lr}$ 
shows the map $X_l\lra Y$ is right orthogonal to $(P)^{l}$, i.e.~belongs to $(P)^{lr}$ as required. 
An argument of this kind is known as Quillen small object argument and originally was 
used to prove Axiom M2  $(cw)(f)$- and $(c)(wf)$-decomposition of model categories. 

The argument shows that under suitable assumptions that a category has enough
limits and colimits, any morphism $X\lra Z$  decomposes as 
$X \xra{(P^l)} Y \xra {(P)^{lr}} Z$, for any class $(P)$ of morphisms.
Here we  take $(P)$ to be the class of proper morphisms.


\vskip15pt\noindent

We end our discussion of compactness with the following rather vague considerations;
we hope they might suggest the reader something  about the arrow notation (calculus) 
appropriate for topology. We admit that what we say below is very vague.

\subsubsection{Compactness as being uniform. $\forall\exists\implies\exists\forall$}
Often an application of compactness is as follows. 
We know that certain choices can be made for each value of parameters;
if we also know that the parameters vary over a compact domain, then 
we may assume that these choices are uniform, 
i.e.~that they do not depend on the value of the parameters.
Put another way, compactness allows to change the order of quantifiers 
$\forall\exists\implies\exists\forall$ in certain formulas. See Appendix~\ref{app:AEEA}
for a list of examples.\footnote{For a discussion  see Remark~8 of
\href{http://mishap.sdf.org/mints/expressive-power-of-the-lifting-property.pdf}{[Gavrilovich, Lifting Property]} 
}

The subsection of [Stacks Project, I.5\S15, tag 005M] dealing with 
the Bourbaki characterisation of proper maps 
starts with a lemma of this kind:
  \begin{lemma}[Tube lemma]
   \label{lemma-tube}
   Let $X$ and $Y$ be topological spaces.
   Let $A \subset X$ and $B \subset Y$ be quasi-compact subsets.
   Let $A \times B \subset W \subset X \times Y$ with $W$
   open in $X \times Y$. Then there exists opens $A \subset U \subset X$
   and $B \subset V \subset Y$ such that $U \times V \subset W$.
   \end{lemma}
In a somewhat more old-fashioned way,
this lemma can be reformulated as follows:
  \begin{lemma}[Tube lemma]
   \label{lemma-tube}
   Let $X$ and $Y$ be topological spaces.
   Let $A \subset X$ and $B \subset Y$ be quasi-compact subsets.
   Let $A \times B \subset W \subset X \times Y$. 

   If for each pair of points $a\in A$ and $b\in B$ 
   we can pick neighbourhoods $U\ni a$ and $V\ni b$ such that 
   $U\times V\subset W$, then we can do so uniformly in $a\in A$ and $b\in B$, 
   i.e.~such that $U$ and $V$ do not depend on $a$ and $b$. 
   
   As a formula, this could be expressed as change of order of quantifiers:
   $$\frac{
\forall a\in A \forall b \in B\, \exists U\subset X \, \exists V \subset Y \,(U\times V\subset W\text{ and }a\in U\text{ is open and }b\in V\text{ is open})    
}{
 \exists U\subset X\, \exists V\subset Y\, \forall a\in A \forall b \in B \,(U\times V\subset W\text{ and }a\in U\text{ is open and }b\in V\text{ is open})
}$$
 \end{lemma}

The following example of change or order of quantifiers is simpler but perhaps more telling.

For a connected topological space $X$, the following are equivalent:
\bi
\item Each real-valued function on $X$ is bounded
\item $\forall x \in K \exists M  ( f(x) < M ) \implies \exists M \forall x \in K  ( f(x) < M ) $
\item $\emptyset \longrightarrow  K \rtt   \sqcup_{ n\in\NN} (-n,n) \longrightarrow   \RR$\\
here  $\cup_ n (-n,n) \longrightarrow   \RR$ denotes the map to the real line
from the disjoint union of intervals $(-n,n)$ which cover it.
Note this is a standard example of an open covering of $\RR$ which
shows it is not compact.
\ei

The following is even more vague. 

\subsubsection{ 
``An open covering has a finite subcovering'' 
} 

Mathematically, this reformulation is based on the following observation:

\def\ooo{\infty}
\begin{quote}
a space $K$ is compact iff for each open covering $U$ of $K$, 
the subset $K$ is closed in $K\cup\{\ooo\}$ in the topology generated 
 elements of $U$ as {\em closed} subsets. 
\end{quote}

This lets us express being {\em finite} with the help of 
the notion of the topology generated by a family of sets.


[Hausdorff, Set theory] denotes by $U(x)$ a neighbourhood of a point $x$, 
which suggests viewing $U(x)$ as a (possibly multivalued) function of a point $x$ ; 
We'd like to develop ``arrow'' notation where this would be expressed as
$$\{x\}\longrightarrow  K\xra{(U(x))}  \{x\rightarrow y\}\ \ \ \ \ \ \ \ \ (*)$$
here it is implicit that $x$ maps to $x$ by the composition of the two arrows;
``$(x)$'' in ``$U(x)$'' signifies that $U(x)$ depends on $x$.

Changing a single symbol ``$\rightarrow$'' into ``$\leftarrow$'' leads us to consider 
elements of $U$ as closed subsets of $K$:
$$\{x\}\longrightarrow  K\xra{(U(x))}  \{x\leftarrow y\}\ \ \ \  (**)$$

We'd like to assume (or require) that $(**)$ inherits some properties of $(*)$, 
in the arrow calculus we'd like to define; 
this would be what corresponds to considering {\em the topology generated by}.
 



\subsubsection{Summary.} These three examples suggest 
that orthogonality, or $\rtt$-negation,  
has a surprising generative power 
as a means of
defining natural elementary mathematical concepts. 
In  Appendix~\ref{app:rtt-examples} and Appendix~\ref{app:rtt-top} we give a number of examples in
various categories, in particular showing that many standard elementary notions
of abstract topology can be defined by applying the lifting property to simple
morphisms of finite topological spaces. 
Examples in topology include the notions of: compact, discrete, connected, and
totally disconnected spaces, dense image, induced topology, and separation axioms.
Examples in algebra include: finite groups being 
nilpotent, solvable, torsion-free, $p$-groups, and prime-to-$p$ groups;
injective and projective
 modules; injective, surjective,
and split homomorphisms.

\subsection{Hausdorff axioms of topology as diagram chasing computations with finite categories
}
We shall now reformulate the axioms of a topology in a form 
almost ready to be implemented in a theorem prover based on diagram chasing. 


Early works talk of topology in terms of {\em neighbourhood} systems $U_x$
where $U_x$ varies though 
 {\em neighbourhoods of points} of a topological space. This is how
the notion of topology was defined by Hausdorff; in words of [Bourbaki], 
``We shall say that a set $E$ carries a topological structure whenever we have 
associated with each element of $E$, by some means or other, a family 
of subsets of $E$ which are called neighbourhoods of this element - provided 
of course that these neighbourhoods satisfy certain conditions (the axioms 
of topological structures).'' 

A neighbourhood $U_x$ of a point $x\in E$ determines two functions
$$
\{x\}\longrightarrow X\xrightarrow{\,\,\, U\,\,\,}\{x{\small\searrow}x'\}$$


This simple observation allows us to show that 
the axioms of topology formulated in the more modern
language of open subsets 
can be seen as diagram chasing rules for manipulating diagrams
involving notation such as
$$   \{x\}\longrightarrow X\ \ \ \ X\longrightarrow \{x{\small\searrow}y\}\ \ \ \  X\longrightarrow \{x\leftrightarrow y\}   $$
in the following straightforward way. 

As is standard in category theory, identify a point $x$ of a topological space $X$
with the arrow $\{x\}\longrightarrow X$, a subset $Z$ of $X$ with the arrow $X\longrightarrow \{z\leftrightarrow z'\}$,
and an open subset $U$ of $X$ with the arrow $X\longrightarrow \{u{\small\searrow}u'\}$.
With these identifications, the Hausdorff axioms of a topological space become
rules for manipulating such arrows, as follows.

{\em Both the empty set and the whole of \ensuremath{X} are open} says that the compositions
$$ X\longrightarrow \{c\}\longrightarrow \{o{\small\searrow}c\}\text{ and }X\longrightarrow \{o\}\longrightarrow \{o{\small\searrow}c\}  $$
behave as expected (the preimage of \{o\} is empty under the first map,
and is the whole of \ensuremath{X} under the second map).

{\em The intersection of two open subsets is open} means the arrow
$$   X\longrightarrow \{o{\small\searrow}c\}\times\{o'{\small\searrow}c'\} $$
 behaves as expected (the ``two open subsets'' are the preimages of points $o\in\{o{\small\searrow}c\}$ and $o'\in\{o'{\small\searrow}c'\}$;
``the intersection'' is the preimage of $(o,o')$ in
$\{o{\small\searrow}c\}\times\{o'{\small\searrow}c'\}$ ).

Finally, {\em a subset $U$ of $X$ is open iff each point $u$ of $U$ has an open
neighbourhood inside of $U$}
 corresponds to the following diagram chasing rule:
\vskip4pt

for each arrow $X\xra{\xi_U}  \{U\leftrightarrow \bar U \}$ it holds\\
 $ \xymatrix{ {\ \ \ \ \ \ \ \ }  & \{U\ra\bar U\} \ar[d] \\
{\ \ \ \ X\ } \ar[r]|{\xi_U} \ar@{-->}[ur] & { \{U\leftrightarrow \bar U \} }}$
\ \ \ \ iff  for each $\{u\}\lra X$, \ \ \ \ \ \ \
 $ \xymatrix{ {\{u\}} \ar[r]  \ar[d] & {\{u \ra U \leftrightarrow \bar U\}} \ar[d] \\
{\ \ \ X\ \ \ } \ar[r]|{\xi_U} \ar@{-->}[ur] & { \{u\!=\!U\!\leftrightarrow\! \bar U \} }}$

%
%
%
%
%
 {\em The preimage of an open set is open} says the composition $$
 X\longrightarrow Y\longrightarrow \{u{\small\searrow}u'\}\longrightarrow \{u\leftrightarrow u'\} $$
is well-defined.

We hope that this reinterpretation may help clarify the nature of the axioms of a
topological space, in particular it offers a constructive approach and a diagram chasing 
formalisation of certain elementary arguments, may clarify
to what extent set-theoretic language is necessary, and perhaps help to suggest
an approach to ''tame topology'' of Grothendieck,
i.e.~a foundation of topology "without false
problems" and "wild phenomena" "at the very beginning".

\section{Topological spaces as simplicial filters.} 
We shall now introduce terminology which we feel allows us to more directly give 
precise meaning to phrases such as ``{\em such and such a property holds for all points sufficiently near $a$}''.
We will do so by ``reading it off'' the informal considerations of [Bourbaki, Introduction] 
of the intuitive notions of limit, continuity and neighbourhood. 

This section may be read independently of the rest of the paper. Unlike the previous section, 
we do not introduce a concise formal syntax to describe the categorical structures that arise.

\subsection{Reading the definition of topology.} 
Now we pretend to directly transcribe the following explanations of Bourbaki of the intuition of topology 
and analysis [Bourbaki, Introduction, p.13]\newline\noindent 
\includegraphics[width=\linewidth]{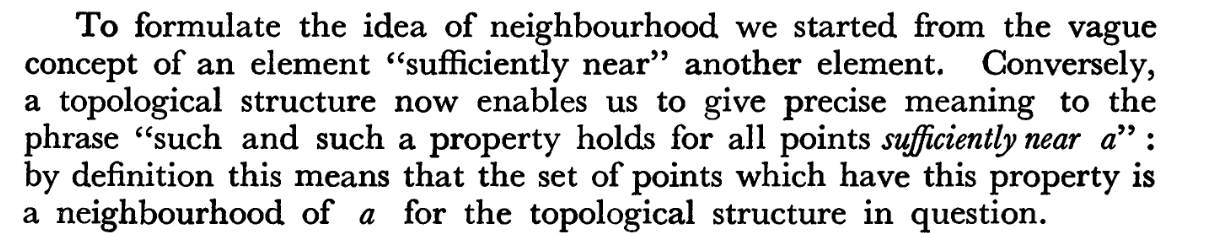}
\includegraphics[width=\linewidth]{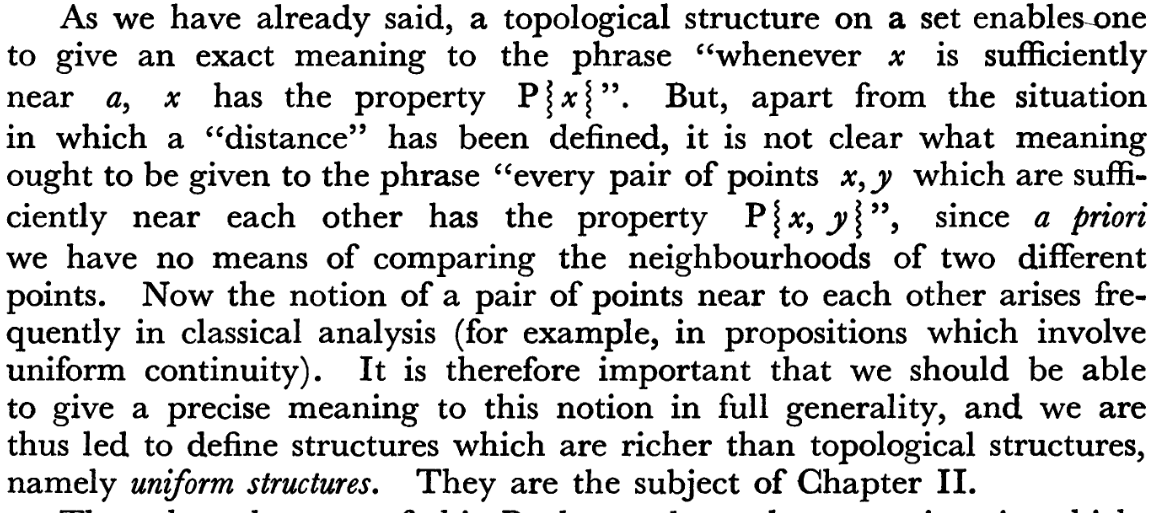}
So let us treat ``whenever $x$ is sufficiently
near $a$, $x$ has the property $P\{x\}$''  as a {\em definition}: 
we {\em define} the exact meaning to the phrase 
``whenever $x$ is sufficiently near $a$, $x$ has the property $P\{x\}$''
to be ``$P$ is topoic'' and introduce the notion of a topoic structure.
We shall say that a set $E$ carries a {\em topoic structure} whenever we have 
defined, by some means or other, a family 
of subsets of $E^n$, $n>0$, which are called {\em topoic subsets of $E^n$} - provided 
that these topoic subsets satisfy certain conditions (the axioms 
of topoic structures).

Similarly,  in the situation in which a ``distance'' has been defined, 
to the phrase ``{\sf every pair of points $x,y$ which are sufficiently near each other has 
the property $P \{x, y\}$}'' we give the precise meaning $P$ is topoic
with respect to the topoic structure associated to the distance (metric) on a set $E$.

A topological argument often relies on consequently choosing ``sufficiently near'' 
points; in this case we expect that 
it implicitly constructs a topoic subset of $E\times .. \times E$.

Sometimes an argument chooses points not consequently, and we hope that often enough 
it implicitly constructs  a topoic subset of $E\times .. \times E$,
albeit in a topoic structure not associated with a topological structure
and possibly specific to the argument.

Evidently the choice of axioms to be imposed 
is to some extent arbitrary, but also depend on whether we consider 
topoic structures associated with topological or metric structures.
In this note we do not discuss this choice.

The conditions on topoic structures associated with a topological structure 
will enable us to define a simplicial object in the category of filters,\footnote
{We find it convenient to allow filters where the empty set is big,
i.e.~we allow the filter of all subsets of a set.

The category $\Filt$ of filters can be thought in three equivalent ways: 
(i) it is a full subcategory of the category of topological spaces whose
objects are spaces such that a subset containing a non-empty open subset is
necessarily open (ii) its objects are sets equipped with a finitely additive
measure taking only two values $0$ and $1$ and such that a subset of a measure
0 set has necessarily measure 0; morphisms are measurable maps preserving the measure 
(iii)
its objects are sets equipped with a collection of subsets called {\em big}
such that the intersection of two big subsets is big and a subset containing a
big subset is necessarily big as well; morphisms are maps such that the
preimage of a big subset is necessarily big.

(ii) suggests that one may also consider the category $\FFilt$ of filters
localised as follows: two maps $f,g:X\lra Y$ are considered equal as morphisms iff they are
equal almost everywhere, i.e.~the subset $\{x: f(x)=g(x)\}$ is big in $X$.}
i.e.~a functor 
$$\ttt(X): \Ord{\omega}^{op} \lra \Filt$$
where  
$\Filt$ is the category of filters, and
 $\Ord{\omega}$ denotes the category of  categories corresponding to finite linear orders
$$\bullet_1\lra .. \lra \bullet_n,\ 0\leq n <\omega.$$

The conditions on topoic structures associated with a metric structure 
will enable us to define in \S\ref{met:filt} a simplicial object in the category of filters,
a functor 
$$\mU(M): \Ord{\omega}^{op} \lra  \Filt$$ 
which factors as
$$\Ord{\omega}^{op} \lra FinSets^{op} \lra \Filt$$ 
where $FinSets^{op}$ is the category of finite sets. 

Thereby we shall obtain fully faithful embeddings of the category of metric spaces (uniform spaces)
and uniformly continuous maps 
and that of topological spaces 
in the category of simplicial filters
$$\ttt: Top \subset Func( \Ord{\omega}^{op} , \Filt)$$
$$\mU: UniformSpaces \subset Func( \Ord{\omega}^{op} , \Filt)$$

Let us now define a topoic structure on a set $E$ 
associated with a topological structure on $E$. 

\subsubsection{\label{def:topoic-top}Topoic structure of a topological space.}
Let $X$ be a topological space. Call a property (subset) $P\subseteq X\times X$ {\em topoic} iff 
$(a,x)\in P$ holds whenever  $x$ is sufficiently near $a$, i.e. 
for each point $a\in X$ there is a neighbourhood $U_a$ such that $(a,x)\in P$ whenever $x\in U_a$.
Call a property $P\subseteq X^n$ {\em topoic} iff
we can ensure that  $(x_1,..,x_n)\in P$ provided
we pick 
$x_2$  sufficiently near $x_1$, then pick $x_{3}$  sufficiently near $x_2$, then
... then pick $x_{n}$  sufficiently near $x_{n-1}$, given any $x_1\in X$, i.e. \\
\begin{small}
for each point $x_1\in X$ there is an open neighbourhood $U_{x_1}\ni x_1$ such that\\
for each point $x_2\in U_{x_1}$ there is an open neighbourhood $U_{x_1,x_2}\ni x_2$ such that\\
for each point $x_3\in U_{x_1,x_2}$ there is an open neighbourhood $U_{x_1,x_2,x_3}\ni x_3$ such that\\
for each point $x_4\in U_{x_1,x_2,x_3}$ ....\\
.....\\
for each point $x_n\in U_{x_1,x_2,...,x_{n-1}}$ there is a neighbourhood $ U_{x_1,x_2,...,x_{n-1}}\ni x_n$ such that\\
$(x_1,...,x_n)\in P$. \\
\end{small}
Topoic subsets form a filter (as well as a topology) on $X^n$: 
 $P'\supset P$, $P$ topoic implies $P'$ is topoic, and the intersection of finitely many topoic sets is topoic. 

As noted above, the {\em filter of topoic subsets} allows us to directly speak about ``sufficiently near'' points.
If a topological argument relies on consequently choosing of ``sufficiently near'' points, then we expect that 
it implicitly constructs a topoic subset of $X\times .. \times X$.

\subsubsection{\label{def:topoic-metric}Topoic structure of a metric space.}
Let $M$ be a metric space. Let us define the topoic structure associated 
with metric space $M$: a subset $P\subset M^n$ is topoic iff there is
$\varepsilon>0$ such that $(x_1,..,x_n)\in M$ provided
$dist(x_i,x_j)<\varepsilon$ for each $1\leq i\leq j\leq n$.  Thereby we give
the phrase ``{\sf every pair of points $x,y$ which are sufficiently near each other
has the property $P \{x, y\}$}'' the precise meaning that $P$ is topoic with respect to
the topoic structure associated with the metric (distance) on $M$. 


\subsubsection{Continuity in topological spaces.} Let us now see that the intuitive explanation of
continuity of a function by  [Bourbaki, Introduction, p.13]  transcribes
directly to the language of topoic subsets.
\newline\noindent\includegraphics[width=\linewidth]{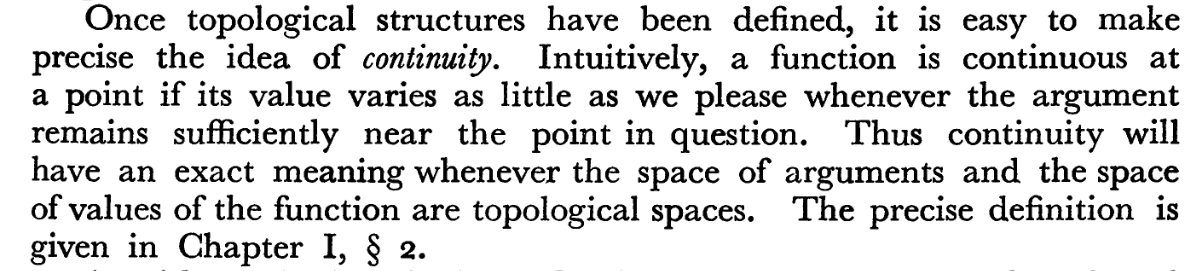}\newline\noindent
%
This reads as: given a topoic subset $W\subset Y\times Y$, we can find a topoic subset 
of $V\subset X\times X$ such that $(f(x_0),f(x))\in W$ provided $(x_0,x)\in V$. 
Here ``given'' corresponds to {\sf as we please}; ``we can find'' to {\sf sufficiently};
the subset $W$ being ``topoic'' corresponds to 
{\sf its value varies as little}; 
and {\sf whenever the argument remains sufficiently {\em near the point in question}} 
corresponds to finding a topoic subset $W$  of $X\times X$ such that  $(f(x_0),f(x))\in W$ provided $(x_0,x)\in V$. 

That is, the map $f\times f:X\times X\lra Y\times Y$ is continuous wrt the topoic filters. 
\subsubsection{Axioms of topology.} Bourbaki reformulate the axioms of topology 
as  Axioms $(V_I-V_{IV})$  stated 
in terms of neighbourhood filters [Bourbaki,I\S1.2], also cf. [ibid, Proposition
2].  Note that the notion of a neighbourhood is all that is need to define
topoic subsets, and let us now try to understand these axioms in terms of
topoic subsets and coordinate maps between $X^n$. 
\newline\noindent\includegraphics[width=\linewidth]{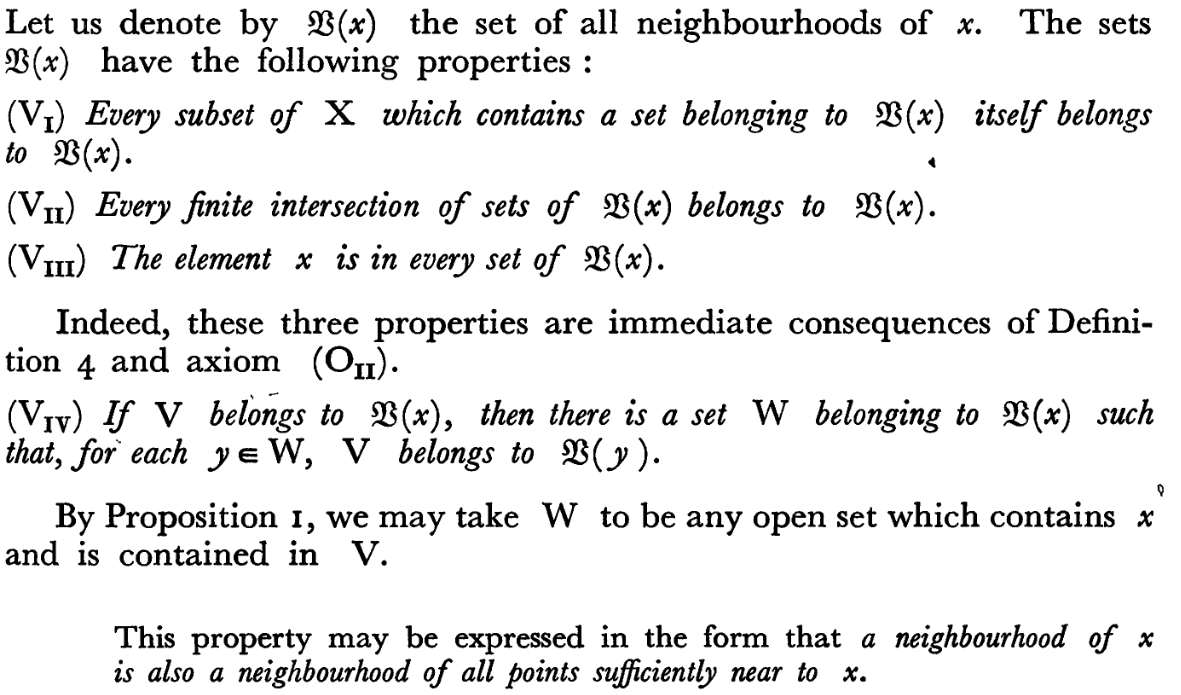}
Axioms $(V_I)$ and $(V_{II})$
say that topoic subsets (as defined above) do indeed form a filter.  
\newline\noindent\includegraphics[width=\linewidth]{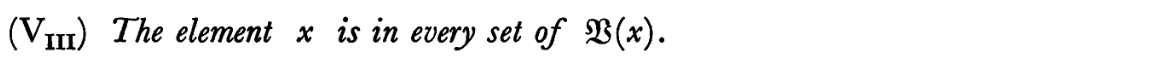}
The the filter on $X$ is antidiscrete, and 
this implies that the diagonal embedding 
$$X\lra X\times X$$ is continuous wrt topoic filters.
\newline\noindent\includegraphics[width=\linewidth]{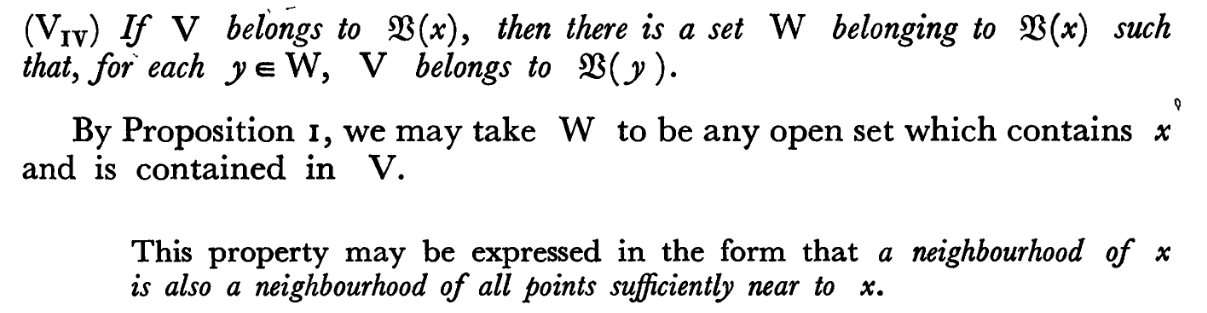}
This means that the map $$X\times X\times X\lra X\times X,\ \ (x_1,x_2,x_3)\mapsto (x_1,x_3)$$
is continuous; to see this, consider the preimage of $\{x\}\times V\subset X\times X$, or rather
a topoic subset $\{x\}\times V \cup (X\setminus\{x\}\times X)\subset X\times X$.
It has to be topoic, and by definition that means that 
for each point $x\in X$, 
 there is a neighbourhood $W$ of $x$ 
such  that, for
 each $y\in W$, $V$ is a neighbourhood of $y$. 
(We may take $W$ to be any open set which contains $x_1$ 
and is contained in $V$.)

 A further verification shows that the maps 
$$X^n\lra X^m,\ (x_1,...,x_n)\mapsto (x_{i_1},...,x_{i_m}),\ 1\leq i_1 \leq ... \leq i_m\leq n$$
have the property that the preimage of a topoic set is topoic. 
\subsubsection{Summing up: a simplicial object of a topological space.}
The coordinate maps between Cartesian powers remind us of the simplicial object of Cartesian powers,
and we are tempted to understand the topological structure as a construction 
of a simplicial object.
And indeed, considerations above show that we obtain a functor
$$\ttt(X): \Ord{\omega}^{op} \lra \Filt$$
where
$\Filt$ is the category of filters, and 
 $\Ord{\omega}$ denotes the category of  categories corresponding to finite linear orders 
$$\bullet_1\lra .. \lra \bullet_n,\ 0\leq n <\omega.$$ 
For a space $X$, the simplicial object $\ttt(X)$ is the object
$$(X_\ttt,X\times X_\ttt, X\times X\times X_\ttt, ...) $$
consisting of Cartesian powers of the set of points of $X$
equipped with the filter of topoic subsets corresponding to the topological structure on $X$
defined in \S\ref{def:topoic-top}.

Continuous maps $f:X\lra X'$ are in one-to-one correspondence with 
natural transformations $\ttt(X)\implies \ttt(X')$, and in fact 
there is a fully faithful embedding of the category of topological spaces
in the category of simplicial filters
$$Top\subset Func( \Ord{\omega}^{op} , \Filt)$$

\subsubsection{\label{met:filt}Metric spaces.}
Consider the topoic structure associated with a metric space $M$.
A straightforward verification shows 
that permutations of coordinates  $M^n\lra M^m, (x_1,...,x_n)\mapsto (x_{i_1},...,x_{i_m})$, $1\leq i_1, ..., \leq i_m\leq n$
have the property that the preimage of a topoic set is topoic, and 
hence we obtain a functor
$$\mU(M): \Ord{\omega}^{op} \lra  \Filt$$ 
which factors as
$$\Ord{\omega}^{op} \lra FinSets^{op} \lra \Filt$$ 
where $FinSets^{op}$ is the category of finite sets. 
This functor sends $n$ to the set $M^n$ equipped with the filter of topoic subsets, i.e.~the filter of subsets containing an $\varepsilon$-neighbourhood of the diagonal, for some $\varepsilon>0$.

Given a mapping $f:M\lra M'$ of sets of points, the condition that 
the preimage of a topoic subset of $M\times M$ is necessarily a  topoic 
subset of $M'\times M'$, 
says that for each $\delta>0$ there is $\varepsilon>0$ such that 
$\dist(f(x),f(y))<\delta$ whenever $\dist(x,y)<\varepsilon$, i.e.
the mapping $f$ is uniformly continuous.

In fact, as is easy to see, this construction also works for uniform spaces,
and we obtain a fully faithful embedding of the category of uniform spaces 
in the category of simplicial filters\footnote{
For more details see 
\href{http://mishap.sdf.org/mints/mints_simplicial_filters.pdf}
{[Gavrilovich, Simplicial Filters]}, in particular Claim~2 which characterises the category of functors
corresponding to uniform spaces.
}
$$\mU: UniformSpaces \subset Func( \Ord{\omega}^{op} , \Filt)$$ 

For a filter $\mathfrak F$, let $\iembE(\FFF)=Hom(n,\FFF)$ denote the simplicial filter
$(\FFF, \FFF\times \FFF, \FFF\times \FFF\times \FFF, ...)$ 
consisting of Cartesian powers of $\FF$ and coordinate maps.

A {\em Cauchy filter} $\FFF$ on a metric space $M$ (cf.~[Bourbaki,II\S3.1,Def.2])
is a filter on the set of points of $M$ such that 
the obvious map $  \iembE(\FFF)  \lra  \mU (M)$ is well-defined.

A {\em Cauchy sequence} in $M$ is a map $\iemb(\NN_{cofinite})\lra \mU(M)$ where
$\NN_{cofinite}$ is the set of natural numbers equipped with cofinite topology (i.e. a subset is closed iff it is finite).

This allows to define various notions of equicontinuity of sequences of functions.

 Let ${X}$ be a topological space, let ${M}$ be a metric space, and
let ${(f_i)_{i \in \NN}}$ be a family of functions $f_i:X\lra M$.

The family ${f_i}$ is {\em equicontinuous} if either of the following equivalent conditions holds:
\bi
\item for every ${x \in X}$ and ${\epsilon > 0}$,
there exists a neighbourhood ${U}$ of ${x}$ such that
${d_Y(f_i(x'), f_i(x)) \leq \epsilon}$ for all ${i \in \NN}$ and ${x' \in U}$
\item the map $\ttt(X)\times \iemb(\{\NN\})\lra \mU(M),\ (x,i)\longmapsto f_i(x)$ is well-defined
\item the map $\ttt(X)\times \iemb(\NN_{cofinite})\lra \mU(M),\ (x,i)\longmapsto f_i(x)$ is well-defined
\ei

    If ${X = (X,d_X)}$ is also a metric space, we say that the family ${f_i}$
is {\em uniformly equicontinuous}
iff either of the following equivalent conditions holds:
\bi \item for every ${\epsilon > 0}$ there exists a ${\delta > 0}$ such that
${d_Y(f_i(x'), f_i(x)) \leq \epsilon}$ for all ${i \in \NN}$ and ${x', x \in x}$ with ${d_X(x,x') \leq \delta}$
\item the map $\mU(X)\times \iemb(\{\NN\})\lra \mU(M),\ (x,i)\longmapsto f_i(x)$ is well-defined
\item the map $\mU(X)\times \iemb(\NN_{cofinite})\lra \mU(M),\ (x,i)\longmapsto f_i(x)$ is well-defined
\ei

The family is {\em uniformly Cauchy} iff
either of the following equivalent conditions holds:
\bi \item for every ${\epsilon > 0}$ there exists a ${\delta > 0}$ and $N>0$ such that
${d_Y(f_i(x'), f_j(x)) \leq \epsilon}$ for all $i,j>N$ and ${x', x \in x}$ with ${d_X(x,x') \leq \delta}$.
\item the map $\mU(X)\times \iembE(\NN_{cofinite}) \lra \mU(M),\ (x,i)\longmapsto f_i(x)$ is well-defined
\ei
Here $\{\NN\}$ denotes the trivial filter on $\NN$ with a unique big subset $\NN$ itself,
and $\NN_{cofinite}$ denotes the filter of cofinite subsets of $\NN$.

This suggest we might reformulate Arzela-Ascoli theorem as something about inner Hom in $\sFilt$,
see Question~\ref{q:ascoli}.

\subsection{Limits as maps to shifted (d\'ecalage) topological spaces.\label{exp:limits}}
We now try to transcribe the explanation of the notion of filter in 
[Bourbaki,I,Introduction].
\newline\noindent\includegraphics[width=\linewidth]{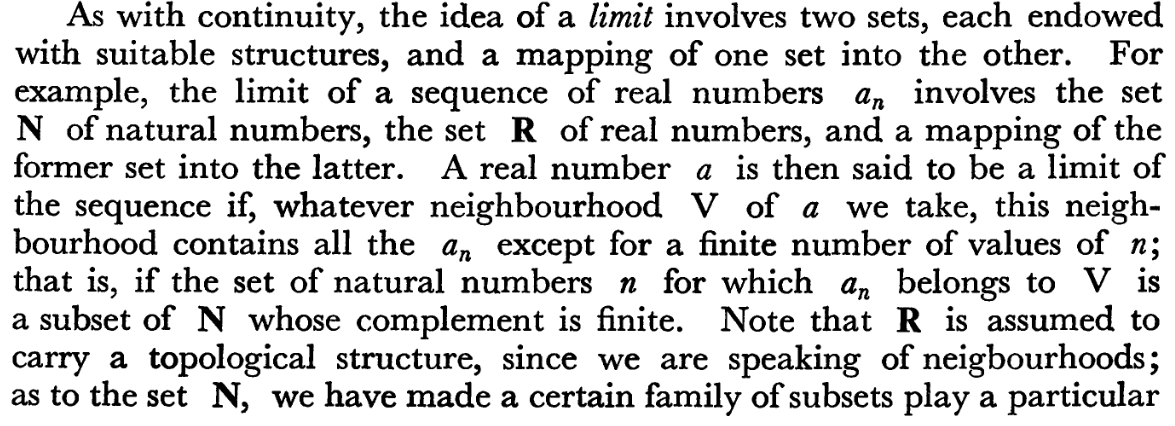}\newline\noindent
 \includegraphics[width=\linewidth]{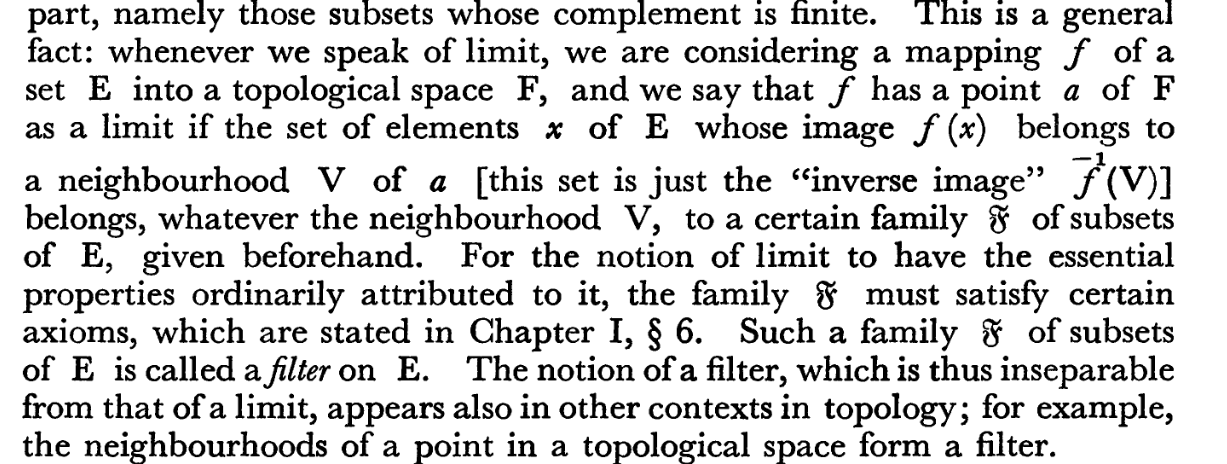}\newline\noindent

So we consider a  mapping $f:E\lra F$ between two sets,  {\sf each endowed with suitable structures}
described rather explicitly: on $F$, it consists of {\sf neighbourhood[s] V of $a$} where $a\in F$ is some fixed point. 
The {\sf suitable structure} on $E$ consists of {\sf a certain family of
subsets, given beforehand} which forms a filter, which we will denote $\FFF$. 
The condition on $f$ says that {\sf ``the inverse image'' $f^{-1}(V)$ belongs
to} the filter, {\sf  whatever the neighbourhood $V$} of $a\in F$.
That is, the mapping $f:E\lra F$ is continuous with respect to the neighbourhood
filter $\mathfrak B(a)$ of point $a\in F$ on $F$
and the given filter $\FFF$ on $E$.

We'd like to have a simplicial map, and indeed this extends to  
$$
 \xymatrix{  {}  & {E_{\mathfrak F}} \ar[r]|{\text{id}} \ar@{->}[d]|{x\mapsto(a,f(x))} & {E_{\mathfrak F}}  \ar[r]|{\text{id}}\ar[l] \ar[d]|{x\mapsto(a,f(x),f(x))} &  {E_{\mathfrak F}}  \ar[r]|{\text{id}}\ar[l] \ar[d]|{x\mapsto(a,f(x),f(x),f(x))} & {} \ar[l] \\ 
{}   &{ {\{a\}\times F}_{\ttt} } \ar[r]\ar[d] & { {\{a\}\times F\times F}_\ttt }  \ar[r]\ar[l]\ar[d] & { {\{a\}\times F\times F\times F}_\ttt }  \ar[r]\ar[l]\ar[d] & {} \ar[l] \\ 
{F}  \ar[r] &{ {F \times F}_{\ttt} } \ar[r]\ar[l] & { {F \times F\times F}_\ttt }  \ar[r]\ar[l] & { {F \times F\times F\times F}_\ttt }  \ar[r]\ar[l] & {} \ar[l] } 
$$ 
Here subscript ${}_\ttt$ means that we consider endowed with the filter of
topoic subsets associated with the topological space $F$.

Let us now introduce notation to describe this diagram. 

For a filter $\FFF$, let $\iemb(\mathfrak F)$ denote the constant simplicial object 
defined by $\iemb(\FFF)(n)=\FFF$ and $\iemb(\FFF)(f)=\id$ for any morphism $f$.

Let $[0]:\Ord{\omega}^{op}\lra \Ord{\omega}^{op}$ denote the endofunctor 
appending a new least element to each linear order, i.e.
$$1< ... < n \,\, \longmapsto\,\, 0 < 1 <...<n$$ 
$$i_1\leq ... \leq i_n \,\, \longmapsto\,\, 0\leq i_1 \leq ... \leq i_n$$ 
Then a topological space $F:\Ord{\omega}^{op}  \lra \Filt$ defines another simplicial object
$F\circ [0]: \Ord{\omega}^{op} \lra \Filt$ which is the bottom row of the diagram.
The middle row we shall denote by  $F\circ[0]_a$; the filter on $\{a\}\times F^n\subset F^{n+1}$ is induced 
by the filter of topoic subsets of $F^{n+1}$. 

With this notation, the diagram above gives a map $$\iemb(\mathfrak F)\xra{x\mapsto(a,x,x,..)} F\circ[0]_a\lra F\circ[0].$$

Let us ponder further; notation $(a,x,x,...)$ is implicitly infinite and thus
somewhat unsatisfactory for us.
 
Setwise, i.e.~if we forget filters,  $F\circ[0]_a$ is an
object of Cartesian powers and coordinate maps, and such objects form a full
subcategory which appears important. Hence, we are tempted to think of the
limit as a decomposition of a map from an ``identity'' simplicial object,
i.e.~an object of form $\iemb(\FFF)$, to the shifted topological space
$F\circ[0]$ via the subcategory consisting of objects which are {\em setwise}
objects of Cartesian powers and coordinate maps. 

And indeed, the following remark shows that such a decomposition determines a point $a\in X$ and is almost of the form above.

\begin{remark} 
Let $\mathcal X$ and $\mathcal Y$ be objects of $\sFilt$ such that are setwise objects of Cartesian powers of sets $X$ and $Y$ and coordinate maps.
Then a morphism $\mathcal X\circ [0]\lra \mathcal Y$ determines a point of $a\in Y$ and is necessarily of form 
$(x_1,..,x_n)\mapsto (a,f(x_1),..,f(x_n))$ for some function $f:X\lra Y$. To see this, take any two points $x_1,x_2\in X$ 
and consider the coordinate projections $(a,y_1)=f(x_1)$ and $(a,y_2)=f(x_2)$
of  $(a,y_1,y_2)=f(x_1,x_2)\in Y\times Y\times Y$ on the second and third
coordinate.

This suggests that the limit of a filter $\FFF$ on a space $X$ might be viewed as something like
factorisaton of $\iemb(\FFF)\lra \ttt(X)\circ [0]$ via the subcategory of objects of Cartesian powers.  
\end{remark}

\begin{remark}\label{ex:paths_as_Aut} Note that one may identify the real interval $[0,1]$ with the set of endofunctors 
$\Ord{\omega}\xra{\,\tau\,} \Ord{\omega}$ such that $\tau(n)=n+1$ for any $n\in Ob  \Ord{\omega}$. To see this, note that such an endofunctor 
is given by inserting a new element in each finite linear order in a compatible manner. 
This temps us to think of the space of paths. 
\end{remark}

\subsection{\label{cwf-decompositions}Path and cylinder spaces and Axiom M2(cw)(f) and M2(c)(wf) of Quillen model categories}

In the category of topological spaces, there is a simple but very useful
way to turn an arbitrary map into either a fibration or cofibration. 
It it captured by Axiom M2 of model categories which requires that 
each map decomposes as a composition of a cofibration and a fibration, and any one of them
may also be required to be a weak equivalence.  

In notation,\\ 
$\xymatrix{ 
& {} &  {A\times [0,1] \sqcup_A B 
} \ar@{<--}[ld]|{(c)}\ar@{-->}[rd]|{(wf)} & {} \\ 
&A\ar[rr] & {} & B}$
$\xymatrix{ 
& {} & A\times_B B^I \ar@{<--}[ld]|{(cw)}\ar@{-->}[rd]|{(f)} & {} \\ 
&A\ar[rr] & {} & B}$

Figure 3 gives drawings
 representing these decompositions in the category of topological spaces.
\newline\noindent\includegraphics[width=\linewidth]{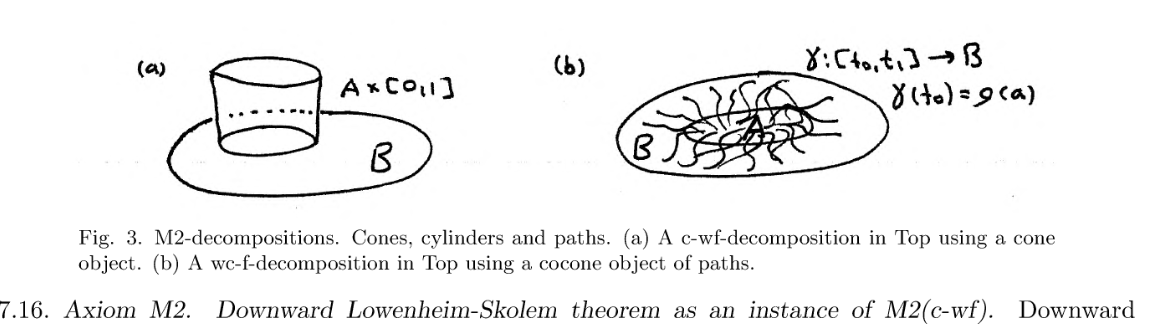}
\subsubsection{($cw$)($f$)-decomposition.}
Let us analyse Figure~3b (cw-f decomposition). 

Recall that, to translate, we care both about intuition 
and algebraic manipulations. 

The construction uses he following algebraic manipulations: 
it considers pairs $(x, \gamma^{const}_{f(x)})$ and 
$(x,\gamma(t_1))$. We ignore paths because these are complicated 
``infinitary'' notions we are unable to express in our language, 
and hence all we are left with are  pairs $(x,f(x))$ where $ x\in A$ and 
$(x,y)$ where $x\in A,y\in B$.
This suggests we look at the following decomposition: 
\newline\noindent
$
 \xymatrix{ 
 {X}  \ar[r]\ar[d]|{(x,f(x))} &{ X\times X } \ar[r]\ar[l]\ar[d]|{(x_1,f(x_1),x_2,f(x_2))} & { X\times X\times X }  \ar[r]\ar[l]\ar[d]|{(x_1,f(x_1),x_2,f(x_2),x_3,f(x_3))} & {} \ar[l]  \\
 {X\times Y}  \ar[r]\ar[d]|{y} &{ X\times Y\times X\times Y } \ar[r]\ar[l]\ar[d]|{(y_1,y_2)} & { X\times Y\times X\times Y\times X\times Y }  \ar[r]\ar[l]\ar[d]|{(y_1,y_2,y_3)} & {}  \ar[l]  \\
 { Y}  \ar[r] &{ Y\times Y } \ar[r]\ar[l] & { Y\times Y\times Y }  \ar[r]\ar[l] & {}  \ar[l] } \\
$
 
What is the filter $\mathfrak F_{\text{cw-f}}$ on the elements of the middle row $X\times Y\times X\times Y \times ...$? 
The diagram suggests that we start with the filter corresponding to the product topology on $X\times Y$. 

Let us use the intuition. In model categories,  
weak equivalences are thought of as equivalences and therefore 
$X$ and $(X\times Y)_{\mathfrak F_{\text{cw-f}}}$ should be very similar for purposes we care about.
Geometric intuition suggests that  we only care about an infinitesimal neighbourhood of
$\{(x,id_x):x\in X\}$ and would prefer our paths to be infinitesimally short.

This motivates us to  modify the topoic filter of the product topology on $X\times Y$ by adding as topoic
``infinitesimal neighbourhoods of $X$''  
$$\bigcup\limits_{x\in X,\, f(x)\in U_{f(x)}\text{ a  neighbourhood}}\,\{x\}\times U_{f(x)} \subset X\times Y$$
That is, we define a new filter $\mathfrak F_{\text{cw-f}}^{n}$ on $(X\times Y)^n$ generated by subsets
$U^n\cap W$ where $U$ is of the form above and $W$ is a topoic subset of $(X\times Y)^n$
with respect to the product topology on $X\times Y$.

\begin{remark} Note that when $X=a$ is a point, this decomposition gives us 
 $\{a\}\xra{(cw)}Y\circ[0]_a\xra{(f)} Y$. 
This suggests that we think of $Y\circ[0]_a$ as the space of infinitesimally short paths starting at $a$;
this would correspond to the intuition that an infinitesimally short path is roughly the same as its endpoint.
\end{remark}

\subsubsection{($c$)($wf$)-decomposition.}
Dually, the {\em (c)(wf)}-decomposition of $X\xra f Y$ leads us to consider 
\newline\noindent
$
 \xymatrix{ 
 {X}  \ar[r]\ar[d]|{(x,f(x))} &{ X\times X } \ar[r]\ar[l]\ar[d]|{(x_1,x_2)} & { X\times X\times X }  \ar[r]\ar[l]\ar[d]|{(x_1,x_2,,x_3))} & {} \ar[l]  \\
 {(X\sqcup Y)}  \ar[r]\ar[d]|{y} &{ (X\sqcup Y)\times (X\sqcup Y) } \ar[r]\ar[l]\ar[d]|{(y_1,y_2)} & { (X\sqcup Y)\times (X\sqcup Y)\times (X\sqcup Y) }  \ar[r]\ar[l]\ar[d]|{(y_1,y_2,y_3)} & {}  \ar[l]  \\
 { Y}  \ar[r] &{ Y\times Y } \ar[r]\ar[l] & { Y\times Y\times Y }  \ar[r]\ar[l] & {}  \ar[l] } \\
$

Intuitively, a (infinitesimal) neighbourhood  of a point $y$ contains points $(x,1)$ whenever $f(x)=y$.

This motivates us to  modify the topoic filter of the disjoint union topology on $X\sqcup Y$
by requiring its topoic subsets to satisfy also the following property:
$x\in P$ whenever $f(x)\in P$, $x\in X$.

That is, we define a new filter $\mathfrak F_{\text{c-wf}}^{n}$ on $(X\times Y)^n$ 
which consists of the subsets topoic wrt the disjoint union topology 
which also satisfy the property that 
$(z_1,..,x_i,...,z_n)\in P$ whenever $x_i\in X$ and $(z_1,..,f(x_i),...,z_n)\in P.$

Arguably, one might find an intuition according to which $(X\sqcup Y)_{\mathfrak F_{\text{c-wf}}}$ is similar to $Y_{\ttt}$.

\section{Open questions and directions for research}

\subsection{Research directions}

Our observations suggest the following broad questions and directions for research.

\subsubsection{Category theory implicit in elementary topology}
We'd like to think of our observations as {\em translations of ideas of
Bourbaki on general topology into a language of category theory appropriate to these ideas}, and
that these ideas (but not notation) are implicit in Bourbaki and reflect their
logic (or perhaps the ergologic in the sense of [Gromov]). 

\begin{enonce}{Question}[Category theory and topological ideas and intuition]
\begin{itemize}
\item Translate more of Bourbaki and some intuitive topological arguments into
the language of category theory and diagram chasing. 
\item Understand how this translation works and in what way it is a {\em translation} 
rather than something new. Formulate what does it mean to say that these
category theoretic constructions are implicit in Bourbaki and find evidence
that indeed they are implicitly there.
\item More speculatively, find evidence that these category theoretic
diagram chasing arguments are implicitly present in the topological intuition of a student,
say by finding correlations between errors of intuition and errors of calculation.
\end{itemize}
\end{enonce}

But is this so and what does it actually mean? 

The goal of our analisys is somewhat reminiscent of the goal of
\href{http://wilfridhodges.co.uk/arabic05.pdf}{[Hodges. Ibn Sina on analysis:
1. Proof search. Or: Abstract State Machines as a tool for history of logic]}
where he ``extract[s] from [the text of Ibn Sina's commentary on a a couple of
paragraphs of Aristotle's Prior Analytics] all the essential ingredients of an
Abstract State Machine for [a proof search] algorithm''. We'd like to think 
that we extract from the text of a couple of paragraphs of Bourbaki
all the  essential ingredients of certain category theoretic constructions.

\subsubsection{Formalisation of topology.} Our translation is unsophisticated and is largely based 
on textual coincidences and correlations between the text and allowed category theoretic manipulations. 
Can these coincidences---and the translation---be found by a machine learning algorithm? A hope is that 
category theoretic manipulations are restrictive enough so that a brute force search for correlations
between (long enough sequences of) allowed category theoretic manipulations and 
the text of Bourbaki may produce meaningful results.  
 
Designing such an algorithm would involve designing a derivation system for the category theoretic constructions used. 

Our reformulations of certain notions of topology in terms of orthogonality (negation) are so concise (several bytes) 
that they can be found, or rather listed, by a brute force search. 
This might be a starting point in designing such an algorithm: first design
an algorithm which can work with the reformulations in terms of 
iterated orthogonals (negations) of maps between finite spaces, 
find correlations between these orthogonals and text of Bourbaki, 
and  single out 
the interesting notions obtained by iterated negation (orthogonals) from very simple morphism.
These notions should include quasi-compactness, denseness, connected etc. 

\begin{enonce}{Question}[Category theory and topological ideas and intuition]
\begin{itemize}
\item Write a short program which extracts diagram chasing derivations from
texts on elementary topology, in the spirit of the ideology of ergosystems/ergostructures.
The texts might include [Bourbaki, General Topology] as well as some informal explanations. 

In particular, it should be able to convert verbal definitions of properties defined by orthogonals
into the corresponding orthogonals.
\item Develop a formalisation of topology based on this translation. 
\end{itemize}
\end{enonce}

\subsubsection{Tame topology and foundations of topology} 
 Does our point of view shed light on the tame topology of Grothendieck 
and allows to develop a foundation of topology
``without false problems'' and ``wild phenomena'' ``at the very beginning''? 
 [Esquisse d'un Programme, translation,\S5,p.33] 

Our approach does seem to avoid certain set-theoretic issues and constructions. 
For example, ultrafilters do appear in our reformulation of compactness, but do so 
only in a combinatorial disguise. Remark~\ref{ex:paths_as_Aut} in \S\ref{exp:limits} suggests a way to think about pathspaces 
without real numbers. 

\begin{enonce}{Question}[Tame topology and foundations of topology.] 
\begin{itemize}
\item Develop elementary topology in terms of finite categories (viewed as
finite topological spaces) and labelled commutative diagrams, with an emphasis on
labels (properties) of morphisms defined by iterated orthogonals ($\rtt$-negation).

\item Develop topology in terms of finite categories, labelled commutative diagrams, 
and simplicial filters. Develop a syntax to describe simplicial filters 
as concise as the syntax of iterated $\rtt$-negation
of maps between finite spaces.  

Does this
lead to tame topology of Grothendieck, i.e. a foundation of topology ``without false
problems'' and ``wild phenomena'' ``at the very beginning'' ?
\end{itemize}
\end{enonce}

 Grothendieck suggests that the following needs to be done first:
\begin{quote}
   Among the first theorems one expects in a framework of tame topology as
I perceive it, aside from the comparison theorems, are the statements which
establish, in a suitable sense, the existence and uniqueness of ``the'' tubular
neighbourhood of closed tame subspace in a tame space (say compact to
make things simpler), together with concrete ways of building it (starting
for instance from any tame map $X \lra {\Bbb R}^+$ having $Y$ as its zero set), the
description of its ``boundary'' (although generally it is in no way a manifold
with boundary!) $\partial T$, which has in $T$ a neighbourhood which is isomorphic to
the product of $T$ with a segment, etc. Granted some suitable equisingularity
hypotheses, one expects that $T$ will be endowed, in an essentially unique
way, with the structure of a locally trivial fibration over $Y$, with $\partial T$ as a
subfibration.
 \end{quote}

\begin{enonce}{Question} 
  Write a first year course introducing elementary topology
  and category theory ideas at the same time, 
  based on the observations above
  and the calculus to be developed. 
  Compactness would be explained with help of all the definitions above;
  Tychonoff theorem follows immediately by a diagram chasing argument 
  from the fact that compactness is given by $\rtt$-negation (orthogonal);
  $\AE\longrightarrow  \EA$ definitions would give students some intuition.

  As a first step, write an exposition aimed at students 
 of the separation axioms and Uryhson Lemma in terms
of the lifting properties.\footnote{
See \url{https://ncatlab.org/nlab/show/separation+axioms+in+terms+of+lifting+properties} for a list of reformulations of the separation axioms
in terms of orthogonals.}
\end{enonce}

 Note that the standard proof of Uryhson lemma can be represented as follows: iterate the lifting property defining normal ($T_4$) spaces
$$   \emptyset \longrightarrow {X} \,\rightthreetimes\,  \{x{\swarrow}x'{\searrow}X{\swarrow}y'{\searrow}y\} \longrightarrow  \{x{\swarrow}x'=X=y'{\searrow}y\}
$$
to prove 
 $$ \emptyset \longrightarrow X\rightthreetimes \{x \swarrow x_1 \searrow ... \swarrow x_n \searrow y \} \longrightarrow  \{x \swarrow x_1 = ... = x_n \searrow y \}$$
  Then pass to the infinite limit to construct a map 
 $ X \longrightarrow  \mathbb{R}$.

\subsubsection{Homotopy and model category structure.}

Let $\FFilt$ be the category $\Filt$ of filters localised as follows:
we consider two morphisms  equal iff they coincide on a big subset of the domain,
i.e.~$f,g:X\lra Y$ are considered equal as morphisms in $\FFilt$ iff
the subset $\{x: f(x)=g(x)\}$ is big in $X$.

\begin{enonce}{Question}[Homotopy theory and model category structure on $\sFilt$ or $\sFFilt$.]
Is there an interesting model category structure on $\sFilt$ or $\sFFilt$? 
Does it lead to interesting homotopy theory of uniform spaces?

In \S\ref{cwf-decompositions} we suggest examples of (cw)(f)- and (c)(wf)-decompositions.
Do the corresponding classes of acyclic cofibrations and fibrations
generate a model structure on $\sFilt$ or $\sFFilt$? 

Does either category have interesting objects corresponding to quotients of topological spaces 
by  a group action? 
\end{enonce}

\subsection{Metric spaces, uniform spaces and coarse spaces}

\subsubsection{Uniform structures}
\begin{enonce}{Question}
Rewrite the theory of uniform structures and metric spaces
 in terms of 
the  category $\sFilt$ of simplicial filters. 
In particular, 
\bi\item reformulate the Lebesgue's number lemma,
partition of unity, and the characterisation of paracompactness 
by A.Stone mentioned by [Alexandroff] (cf.~\S\ref{app:AEEA}). 
\ei\end{enonce}

\begin{enonce}{Question}[Arzela-Ascoli]\label{q:ascoli}
\bee
\item Reformulate various notions of equicontinuity and convergence of a family of functions $f_i:X\lra M$
in terms of maps in $\sFilt$
using e.g. $\iemb(\NN_{cofinite})$, $\Ee(\NN_{cofinite})$, $\iemb(\NN_{cofinite}\cup_{\NN_{cofinite}}\{\infty\} )$,
$\ttt(\NN_{cofinite}\cup_{\NN_{cofinite}}\{\infty\}) $, $\Ee(\NN_{cofinite}\cup_{\NN_{cofinite}}\{\infty\} )$,
 $\ttt(\NN_{cofinite})$, $\ttt(X)$, $\mU(X)$, and $\mU(M)$.
\item Reformulate and prove Arzela-Ascoli theorem in terms something like inner $Hom$ in $\sFilt$
and the lifting properties defining precompactness, compactness etc.
\item Define various function spaces in terms of something like inner $Hom$ in $\sFilt$.
\eee
\end{enonce}

\subsubsection{Large scale geometry.}

The category of quasigeodesic metric spaces and large scale Lipschitz maps 
embeds into another category $\sFmilt$ of simplicial filters, with maps of filters defined 
differently: a $\Fmilt$-morphism of filters maps a small subset into a small subset. 


Let $X$ be a metric space. Call a subset $U$ of $X^n$ {\em small} iff the diameters of tuples in $U$ are uniformly bounded,
i.e. there is a $d=d(U)$ such that for each $(u_1,...,u_n)\in U$, $dist(u_i,u_j)\leq d$ for each $1\leq i,j\leq n$;
this defines a filter on $X^n$. 
Note that coordinate maps $X^n\lra X^m$ have the property that the image of a small subset is necessarily small.
Hence this construction defines a functor ${\mathcal X}:FinSets^{op} \lra \Fmilt$. 
A natural transformation $\mathcal X\lra \mathcal Y$ of functors associated with metric spaces $X$ and $Y$, resp.,
corresponds to a map of metric spaces $f:X\lra Y$ such that
for each $d>0$ there is $D>0$ such that $$dist(f(x'),f(x''))<D\text{ whenever }dist(x',x'')<d, x',x'' \in X.$$ 
For $X$ quasi-geodesic, this is the class of large scale Lipschitz maps.

\begin{enonce}{Question}[Large scale geometry]
Rewrite 
 in terms of
the  category $\sFmilt$ of simplicial filters
the theory of metric spaces and uniformly bounded maps
and the theory of coarse structures (cf.~[Bunke, Engel]).
\end{enonce}

\subsubsection{Group theory} In the category of groups, properties defined by
orthogonals (cf.~\S\ref{app:rtt-top}) include groups being
nilpotent, solvable, torsion-free, $p$-groups, and prime-to-$p$ groups.

This suggests it is worthwhile to try to rewrite group theoretic arguments 
in diagram chasing manner, say the proof that nilpotent groups are solvable,
and try to find a semantics for our notation of finite topological spaces 
in the category of groups.

\begin{enonce}{Question}[Group theory]
\bi
\item Calculate iterated $\rtt$-negation (orthogonals) of interesting morphisms in the category of groups
and find interesting properties defined this way.
\item Find a diagram chasing reformulation of the Sylow theorems.  
\item Find a semantics in the category of groups for the notation introduced in \S\ref{app:top-notation}. 
\ei 
\end{enonce}

To reformulate the Sylow theorem, the following characterisation of inner automorphisms may be of help: 
an automorphism $f:G\lra G$ is inner iff either of the following equivalent conditions hold (cf.~[Schupp,Inn]):
\bi \item $f:G\lra G$ extends to an automorphism of $f':H\lra H$, for any $h:G\lra H$, i.e.~$f\circ h=h\circ f'$. 
\item $f:G\lra G$ extends to an automorphism of $f':H\lra H$, for any $h:H\lra G$, i.e.~$f\circ h=h\circ f'$
\ei 
 
To find a semantics, it would help to find a category which contains both groups and finite preorders. 
One candidate is the category $Cats$ of categories where a group $G$ is identified with a category with 
a single object $O$ such that $Aut(O)=G$.
Intuitively, one may think of the category of groups as analogous to the (sub)category of Hausdorff spaces
in the following way: the interesting example are groups (Hausdorff spaces), yet the big ambient category 
contains useful objects (finite topological spaces) which are very unlike the interesting examples we care about,
but are useful to talk about these examples.

\subsection{Open problems.}

Now we would like to formulate several suggestions with specific details.

\subsubsection{Topology.}

\begin{enonce}{Question} 
 Develop a syntax and a derivation calculus based on $\rtt$-negation,  
arrows, labelled arrows and diagrams, finite topological spaces, and simplicial filters.
Develop an intuition for the calculus as well. 
\bi 
  \item[1.] Standard arguments and definitions in elementary topology 
 should be represented by short formal calculations 
    which are both human readable and computer verifiable. 
  \item[2.] In particular, the calculus should express concisely all the three definitions
 of compactness, and prove their equivalence by short formal calculations. 
\ei
\end{enonce}

\begin{enonce}{Question} 
 Does topological intuition (as developed by a first year student) relate to the formal calculus we'd like to develop?
Note that this might be testable by an experiment, namely it might be possible to test whether
mistakes of intuition correspond to mistakes of calculation. This might even be used to develop the calculus.
\end{enonce}

\begin{enonce}{Question} Write a {\em very short} program which would ``invent'' (generate) the (very) basic theory of topology,
possibly using unstructured input such as the text of (Bourbaki, General Topology).
Our examples suggest that iterating right and left $\rtt$-negation up to 5 times 
and restricting size to 3 or 4 is enough to generate, but not single out,  the notions of compactness, connectedness,
a subset, a closed subset, separation axioms, and some implications between them.

What is the length of a shortest such program? To what extent have the axioms of topology to be hardcoded rather than generated?
\end{enonce}

Let us comment on how such a program may look like.

We observed that there is a simple rule which leads to several notions in
topology interesting enough to be introduced in an elementary course. Can this
rule be extended to a very short program which learns elementary topology?

We suggest the following naive approach is worth thinking about.

The program maintains a collection of directed labelled graphs and certain
distinguished subgraphs. Directed graphs represent parts of a category;
distinguished subgraphs represent commutative diagrams.  Labels represent
properties of morphisms.  Further, the program maintains a collection of rules
to manipulate these data, e.g. to add or remove arrows and labels.

The program interacts with a flow of signals, say the text
of [Bourbaki, General Topology, Ch.1], and seeks correlations
between the diagram chasing rules and the flow of signals.
It finds a "correlation" iff certain strings occur nearby
in the signal flow iff they occur nearby in a diagram chasing rule.
To find "what's interesting",  by brute force it searches
for a valid derivation which exhibits such correlations.  To guide the search
and exhibit missing correlations in a derivation under consideration, it may
ask questions: are these two strings related?  Once it finds such a derivation,
the program "uses it for building its own structure".  Labels correspond to
properties of morphisms.  Labels defined by the lifting property play an
important role, often used to exclude counterexamples making a diagram chasing
argument fail.  In [DeMorgan] we analysed the text of the definitions of
surjective and injective maps showing what such a correlation may look like in
a "baby" case.

A related but easier task is to write a theorem prover doing diagram chasing
in a model category. The axioms of a (closed or not) model category
as stated in [Quillen,I.1.1] can be interpreted as rules to manipulate
labelled commutative diagrams in a labelled category. It appears
straightforward how to formulate a logic (proof system) based on these rules
which would allow to express statements like: Given a labelled commutative diagram,
(it is permissible to) add this or that arrow or label.
Moreover, it appears not hard to write a theorem prover for this logic
doing brute force guided search. What is not clear whether this logic is complete
in any sense or whether there are non-trivial inferences of this form to prove.

Writing such a theorem prover is particularly easy when the underlying category
of the model category is a partial order [Gavrilovich, Hasson] and [BaysQuilder] wrote
some code for doing diagram chasing in such a category. However, the latter
problem is particularly severe as well.

The two problems  are related; we hope they help to clarify the notion
of an ergosystem and that of a topological space.

The following are somewhat more concrete questions. 

\begin{enonce}{Question}\bi
\item 
   Prove that a compact Hausdorff space is normal by diagram chasing;
       does it require additional axioms?
       Note that we know how to express the statement entirely in terms of $\rtt$-negation
        and finite topological spaces of small size.
\item 
Formalise the argument in [Fox, 1945] which implies the category of
       topological spaces is not Cartesian closed; does it apply to $\sFilt$?

  Namely, Theorem \ensuremath{3} [ibid.] proves that
          if \ensuremath{X} is separable metrizable space, \ensuremath{R} is the real line, then
           \ensuremath{X} is locally compact
            iff
           there is a topology on $X^\RR$ such that for any space $T$,
           a function $h:X\times T\longrightarrow \RR$ is continuous iff
           the corresponding function $h^*:T\longrightarrow X^\RR$ is continuous
           (where $h(x,t)=h^*(t)(x))$
    
   Note that here we do not know how to express the statement.
\ei \end{enonce}

\begin{enonce}{Question}
Characterise the interval $[0,1]$, a circle $\mathbb S^1$ and, more generally,  spheres $\mathbb S^n$ using their topological characterisations
provided by the Kline sphere charterisation theorem and its analogues. An example of such a characterisation is
that a topological space $X$ is homomorphic to the circle $\mathbb S^1$ iff
$X$ is a connected Hausdorff metrizable space such that $X\setminus \{x,y\}$ is not connected for any two points $x\neq y\in X$
([Hocking, Young. Topology, Thm.2-28,p.55]); another example is that a topological space $X$ is homomorphic to 
the closed interval $[0,1]$ iff 
$X$ is a connected Hausdorff metrizable space such that $X\setminus \{x\}$ is not connected for exactly two points $x\neq y\in X$
([Hocking, Young. Topology, Thm.2-27,p.54]).
\end{enonce}

\subsubsection{$\rtt$-negation, or orthogonality.}
Call a subcategory $\mathcal A$ of $\mathcal B$ {\em $\rtt^s$-full} iff the value of an $s$-orthogonal 
of a class of morhpisms in $\mathcal A$ does not depend  whether it is calculated in $\mathcal A$ or in $\mathcal B$.
i.e.~$(C)^s_{\mathcal A}=(C)^s_{\mathcal B}$ 
for any class $C$ of morphisms of $\mathcal A$, where $s\in \{l,r\}^n$ is a string. 

\begin{enonce}{Question} 
 Calculate left and right $\rtt $-negations and generalisations, e.g. 
   ($C)^{r}$, ($C)^{l}$, ($C)^{rl}$, ($C)^{ll}$, ($C)^{rr}$, ($C)^{llr}$, ... 
   for various simple classes of morphisms in various categories, 
   e.g. morphisms of finite topological spaces or finite groups. 
   
   Develop abstract theory of the lifting property. Find examples of $\rtt^s$-full subcategories.


\end{enonce}

\subsubsection{Compactness as being uniform.}
In \S\ref{app:AEEA} we observe that a number of consequences
   of compactness can be expressed as a change of order of quantifiers in a
   formula,
   i.e. are of form 
   $\forall \exists \phi(...)\implies\exists\forall \phi(...)$
   namely that a real-valued function on a compact is necessarily bounded,
   that a Hausdorff compact is necessarily normal,
   that the image in \ensuremath{X} of a closed subset in $X\times K$ is necessarily closed,
   Lebesgue's number Lemma, and paracompactness.

   Such formulae correspond to inference rules of a special form,
   and we feel a special syntax should be introduced to state
   these rules.

   For example, consider the statement that
   "a real-valued function on a compact domain is necessarily bounded".
   As a first order formula, it is expressed as
$$
   \forall \ensuremath{x} \in \ensuremath{K} \exists M ( f(x) \leq  M ) \Longrightarrow  \exists M \forall \ensuremath{x} \in \ensuremath{K} ( f(x) \leq  M )
$$
   Another way to express it is:
$$
   \exists M:K\longrightarrow \RR \forall \ensuremath{x} \in \ensuremath{K} ( f(x) \leq  M(x) ) \Longrightarrow  \exists M \in \RR \forall \ensuremath{x} \in \ensuremath{K} ( f(x) \leq  M )
$$
   Note that all that happened here is that a function $M:K\longrightarrow \RR$
   become a constant $M \in \RR$, or rather
   expression "M(x)" of type $K\longrightarrow M$ which used (depended upon) variable "x"
   become expression "M" which does not use (depend upon)  variable "x".
   We feel there should be a special syntax which would allow
   to express above as an inference rule {\em removing dependency of "M(x)" on "x"},
   and this syntax should be used to express consequences of compactness
   in a diagram chasing derivation system for elementary topology.

   To summarise, we think that compactness should be formulated
   with help of inference rules for expressly manipulating which variables are 'new', 
   in what order they 'were' introduced, 
   and what variables terms depend on, e.g.~rules replacing a term t(x,y) by term t(x).

   Something like the following:

\begin{verbatim}
   ... f(x) =< M(x) ...
   --------------------
   ...  f(x) =< M   ...
\end{verbatim}

\begin{enonce}{Question} 
 In \S\ref{app:AEEA} we give several examples where consequences of compactness are expressed as 
change of order of quantifiers $\AE\longrightarrow  \EA$.
\bi\item        
Is there a theorem generalising these examples? 
\item Is there a proof system which allows to formulate inference rules corresponding to these reformulations? 
\ei
\end{enonce}  

\subsection{Open problems.}

Here we formulate precise questions one may ask.
The choice of these questions is somewhat arbitrary. 

\begin{enonce}{Question}[Iterated orthogonals in $Top$]
\bi
\item
Are there finitely many different iterated orthogonals of the form 
$\{\emptyset\lra\{\bullet\}\}^s$ where $s\in \{l,r\}^{<\omega}$?

More generally, are there finitely many different classes obtained from 
$\{\emptyset\lra\{\bullet\}$ by repeatedly passing to left or right orthogonal $C^l$ or $C^r$ or 
the subclass $C_{<n}$ of morphisms between spaces of size at most $n$? 

Is there an algorithm which decides whether two such classes are equal? 

\item
Find the shortest expressions (Kolmogoroff complexity) of various topological notions.

\item 
Is $$  
 ((\{ \{o\}\lra \{o\ra c\}\}^{r})_{<5})^{lr}$$
the class of proper maps? 

\item Calculate\footnote{
For motivation see Remark~5 of 
\href{http://mishap.sdf.org/mints/expressive-power-of-the-lifting-property.pdf}{[Gavrilovich, Lifting Property]}
} $$  
((\{c\}\longrightarrow \{o{\small\searrow}c\})^r_{<5})^{lr},\ \ \ \ \  ((\{c\}\longrightarrow \{o{\small\searrow}c\})^{l
rr}$$  
 $$  (\{a{\small\swarrow}U{\small\searrow}x{\small\swarrow}V{\small\searrow}b\}\longrightarrow \{a{\small\swarrow}U=x=V{\small\searrow}b\})^{lr} $$
\ei 
Note the orthogonals  may depend on the category they may are calculated in, 
which is either $Top$ or $\sFilt$.
\end{enonce}


 A number of elementary topological properties can be defined by, in a sense, combinatorial expressions,
   by taking iterated orthogonals in $Top$ 
   of a single 
   morphism between finite topological spaces [Gavrilovich, Lifting Property]. Calculate 
     these expressions in $\sFilt$ using the embedding $\ttt:Top\lra \sFilt$. Note that this would give properties of both topological spaces 
     and metric spaces.  Do they define the same properties of topological spaces? Do they provide an interesting analogy between topological spaces 
and metric spaces, e.g. compactness [Bourbaki, I\S10.2, Thm.1(d), p.101] and completeness [Bourbaki, II\S3.6, Prop.11]? 

The following is an example of a precise conjecture.

\begin{enonce}{Question}[Compactness and completeness]
\bi
\item Calculated in $Top$, is 
 $$((\{ \{o\}\lra \{o\ra c\}\}^{r})_{<5})^{lr}$$
the class of proper maps?
\item 
Is the following true in n $\sFilt$ or $\sFFilt$? 
 \bee
\item A Hausdorff space $X$ is compact iff $$\ttt\left(X \lra \{\bullet\}\right) \in \left((\{\ttt\left(\{o\}\lra \{o\ra c\}\right)\}^{r})_{<5}\right)^{lr} $$
\item A metric space $M$ is complete iff $$\mU\left(M\lra \{\bullet\}\right) \in \left((\{\ttt\left( \{o\}\lra \{o\ra c\}\right)\}^{r})_{<5}\right)^{lr} $$
\eee
\item Does the value of an orthogonal depend whether it is calculated in $Top$ or $\sFilt$? For example, is it true that for any morphism $f$ of  finite topological spaces, 
\bi\item[] $\{f\}^{lr}_{Top} = \{\ttt(f)\}^{lr}_{\sFilt}\cap Top $ and
$\{f\}^{rl}_{Top} = \{\ttt(f)\}^{rl}_{\sFilt}\cap Top $ ?
\ei\ei
\end{enonce}


Recall $\FFilt$ is the category $\Filt$ of filters localised as follows:
we consider two morphisms  equal iff they coincide on a big subset of the domain,
i.e.~$f,g:X\lra Y$ are considered equal as morphisms in $\FFilt$ iff
the subset $\{x: f(x)=g(x)\}$ is big in $X$. 

\begin{enonce}{Question}[Is $\sFilt$ or  $\sFFilt$ a model category?]
Let $(cw)_0$, $(f)_0$, $(c)_0$ and $(wf)_0$ be the classes of maps 
arising in the examples of  M2 $(cw)(f)$- and $(c)(wf)$-decompositions
suggested in \S\ref{cwf-decompositions}. 

Do classes  $(cw)_0^{lr}$, $(f)_0^{rl}$, $(c)_0^{lr}$ and $(wf)_0^{rl}$
define a model structure on $\sFilt$ or on $\sFFilt$? 

Does it induce one of the usual model category structures on the subcategory $Top$ of $\sFilt$?
\end{enonce}

\begin{enonce}{Question}[Orthogonals in group theory]
\bi
\item Is the class of finite CA-groups or CN-groups defined by a natural lifting property, say
 as an iterated orthogonal of a single homomorphism?
   Recall that a group is a CA-group, resp. CN-group, iff the centraliser of a non-identity
   element is necessarily abelian, resp. nilpotent.

\item 
 Calculate iterated left and right $\rtt $-negations and generalisations, e.g. 
   ($C)^{r}$, ($C)^{l}$, ($C)^{rl}$, ($C)^{ll}$, ($C)^{rr}$, ($C)^{llr}$, ... 
   for various simple classes of morphisms in various categories, 
   e.g. morphisms of finite topological spaces or finite groups. 
\item Reformulate the Feit-Thompson odd group theorem as inclusion of orthogonals.\footnote
{For a partial reformulation see~\href{http://mishap.sdf.org/mints/mints-expressing-odd-subgroup-theorem-with-diagrams.pdf}{Gavrilovich, Expressing the statement of the Feit-Thompson theorem with diagrams in the category of finite groups]}.} 
\ei   
\end{enonce}


\section{Appendix}

\subsection{\label{app:sur-and-in}Surjection and injection: an example of translation and orthogonals}

This section is part of a note\footnote{See \href{http://mishap.sdf.org/mints/mints-lifting-property-as-negation-DMG_5_no_4_2014.pdf}{[Gavrilovich, DMG]}. I thank Vladimir Sosnilo for help with the exposition.} written for  {\em The De Morgan Gazette}  to demonstrate that some natural definitions are
lifting properties relative to the simplest counterexample, and to suggest a way to ``extract'' these lifting properties
from the text of the usual definitions and proofs. The exposition is in the form of a story and aims to be self-contained
and accessible to a first year student
who has taken some first lectures in naive set theory, topology, and who has heard a definition of a category.
A more sophisticated reader  may find it more illuminating to recover our formulations herself
from reading either the abstract, or the abstract and the opening sentence of the next two sections.
The displayed formulae and Figure~4(a) defining the lifting property provide complete formulations of our theorems
to such a reader.

\subsubsection{\label{surinj} Surjection and injection}

We try to find some ``algebraic''  notation to (re)write the {\em text} of the
definitions of surjectivity and injectivity of a function, as found in any standard
textbook.
We want something very straightforward and syntactic---notation for what we
(actually) say, for the text we write,  and not for its meaning, for who knows
what meaning is anyway?\\

\bd\item[\textbf{(*)}${}_{\textbf{words}}$] ``A function $f$ from $X$ to $Y$ is {\em surjective} iff for every element $y$ of $Y$ there is an element $x$ of $X
$
such that $f(x)=y$.''\\
\ed
A function from $X$ to $Y$ is an arrow $X\lra Y$. Grothendieck taught us that a point, say ``$x$ of $X$'', is
(better viewed as) as  $\{\bullet\}$-valued point, that is an arrow $$\{\bullet\}\lra X$$
from a (the?) set with a unique element;
similarly ``$y$ of $Y$'' we denote by
an arrow $$\{\bullet\}\lra Y.$$ Finally, make dashed the arrows  required to  ``exist''. We get the diagram Fig.~1(b)
without the upper left corner; there ``$\{\}$'' denotes the empty set with no elements listed inside of the brackets.\\
\bd
\item[\textbf{(**)}${}_{\textbf{words}}$] ``A function $f$ from $X$ to $Y$ is {\em injective} iff  no pair of different points of $X$ is
sent to the same point of $Y$.''\\
\ed
``A function $f$ from $X$ to $Y$'' is an arrow $X\lra Y$.  ``A pair of points'' is a
$\{\bullet,\bullet\}$-valued point, that is an arrow $$\{\bullet,\bullet\}\lra X$$ from a two element set;
we ignore ``different'' for now.
``the same point of $Y$'' is an arrow $\{\bullet\}\lra Y$. Represent ``sent to'' by an arrow 
\[
\{\bullet,\bullet\}\lra \{\bullet\}.
\]
What about ``different''? If the points are not ``different'', then they are ``the same'' point of $X$,
and thus we need to add an arrow representing a single point of $X$,
that is an arrow 
\[
\{\bullet\}\lra X.
\]

Now all these arrows combine nicely into diagram Figure~4(c);  
however, our analysis does not necessarily
makes it clear that the diagonal arrow needs to be denoted differently. 
How do we read it? We want this diagram to have the meaning of the sentence
(**)${}_{\text{words}}$ above, so we interpret such diagrams as follows:\\

\bd\item[$\boldsymbol{(\rtt)}$]
``for every commutative square (of solid arrows) as shown  there is a diagonal (dashed) arrow making the
total diagram commutative'' (see Fig.~$1(a)$).\\
\ed
(recall that ``commutative'' in category theory
means that the composition of the arrows along a directed path depends only on the endpoints of the path)

Property $(\rtt)$ has a name and is in fact quite well-known [Qui]. It is called {\em the lifting property},
or sometimes {\em orthogonality of morphisms},
and is viewed as the property of the two downward arrows; we denote it by $\rtt$.

Now we rewrite (*)${}_{\text{words}}$ and (**)${}_{\text{words}}$ as:
$$\!\!\mathrm{(*)_{\rtt}}\ \ \ \ \ \ \ \ \ \ \ \ \ \ \ \ \ \ \{\} \lra \{\bullet\} \rtt  X\lra Y$$
$$\mathrm{(**)_{\rtt}}\ \ \ \ \ \ \ \ \ \ \ \ \ \ \ \  \{\bullet,\bullet\}\lra \{\bullet\} \rtt  X\lra Y$$

\setcounter{figure}{3}
\def\rrt#1#2#3#4#5#6{\xymatrix{ {#1} \ar[r]|{} \ar@{->}[d]|{#2} & {#4} \ar[d]|{#5} \\ {#3}  \ar[r] \ar@{-->}[ur]^{}& {#6} }}
\begin{figure}
\begin{center}
\large
$ (a)\ \xymatrix{ A \ar[r]^{i} \ar@{->}[d]|f & X \ar[d]|g \\ B \ar[r]|-{j} \ar@{-->}[ur]|{{\tilde j}}& Y }$
$(b)\  \rrt  {\{\}}  {} {\{\bullet\}}  X {\therefore(surj)} Y $
$(c)\  \rrt {\{\bullet,\bullet\}} {} {\{\bullet\}}  X {\therefore(inj)} Y $
$(d)\  \rrt X {\therefore(inj)} {Y} {\{x,y\}} {} {\{x=y\}}$\ 
\end{center}
\caption{\label{fig4}\normalsize
Lifting properties. Dots $\therefore$ indicate free variables and what property of these variables is being defined;
in a diagram chasing calculation, ``$\therefore(surj)$" reads as: 
given a (valid) diagram, add label $(surj)$ to the corresponding arrow.\newline
 (a) The definition of a lifting property $f\rtt g$: for each $i:A\lra X$ and $j:B\lra Y$
making the square commutative, i.e.~$f\circ j=i\circ g$, there is a diagonal arrow $\tilde j:B\lra X$ making the total diagram
$A\xra f B\xra {\tilde j} X\xra g Y, A\xra i X, B\xra j Y$ commutative, i.e.~$f\circ \tilde j=i$ and $\tilde j\circ g=j$.
 (b) $X\lra Y$ is surjective\newline
 (c) $X\lra Y$ is injective; $X\lra Y$ is an epicmorphism if we forget 
that $\{\bullet\}$ denotes a singleton (rather than an arbitrary object
and thus $\{\bullet,\bullet\}\lra\{\bullet\}$ denotes an arbitrary morphism $Z\sqcup Z\xra{(id,id)} Z$)\newline
 (d) $X\lra Y$ is injective, in the category of Sets; $\pi_0(X)\lra\pi_0(Y)$ is injective, 
 when the diagram is interpreted in the category
of topological spaces.
}
\end{figure}

So we rewrote these definitions without any words at all. Our benefits?
The usual little miracles happen:

The notation makes apparent a similarity of (*)${}_{\text{words}}$ and
(**)${}_{\text{words}}$: they are obtained, in the same purely formal way, from
the two of the simplest arrows (maps, morphisms) in the category of Sets. More is true:
it is also apparent that these two arrows are the simplest {\em counterexamples} to the properties,
and this suggests that we think of the lifting property as  a category-theoretic (substitute for) negation.
Note also that a non-trivial (one which is not an non-isomorphism) morphism never has the lifting property
relative to itself, which fits with this interpretation.

Now that we have a formal notation and the little observation above,
we start to play around looking at simple arrows in various categories,
and also at not-so-simple arrows representing standard counterexamples.
You notice a few words from your first course on topology:
{\em 
$(i)$ connected, $(ii)$ the separation axioms $T_0$ and $T_1$, $(iii)$ dense, $(iv)$ induced (pullback) topology}, and $(v)$ {\em Hausdorff 
}
are, respectively,\footnote{The notation is self-explanatory; for the definition see \S\ref{app:top-notation}.}
\bi
\item[(i):] $$ X\lra\{\bullet\}\rtt \{\bullet,\bullet\}\lra\{\bullet\}$$
\item[(ii):] $$\{\bullet \leftrightarrow \star\}\lra\{\bullet=\star\} \rtt X\lra\{\bullet\}$$  and
$$\{\bullet\ra \star\}\lra\{\bullet=\star\} \rtt X\lra\{\bullet\}$$

\item[(iii):]  $$X\lra Y\rtt \{\bullet\}\lra\{\bullet\rightarrow\star\}$$
\item[(iv):] $$X\lra Y\rtt \{\bullet\ra\star\}\lra\{\bullet\}$$
\item[(v):] $$\{\bullet,\bullet'\}\xra{(inj)} X \rtt \{\bullet\leftarrow\star\rightarrow\bullet'\}\lra \{\bullet\}$$
\ei
here
$$\{\bullet\ra\star\},\; \{\bullet\leftrightarrow\star\},\; \dots\ $$
denote finite preorders, or, equivalently,
finite categories with at most one arrow between any two objects, or finite topological spaces
on their elements or objects,
where a subset is closed iff it is downward closed (that is, together with each element, it contains all the smaller elements).
Thus $$\{\bullet\ra\star\},\; \{\bullet\leftrightarrow\star\} \;\mbox{ and }\;  \{\bullet\leftarrow\star\rightarrow\bullet'\}$$
denote the connected spaces with only one open point $\bullet$,  with no open points, and with two open points
$\bullet,\bullet'$ and a closed point $\star$.
Line (v) is to be interpreted somewhat differently: we consider {\em all}
the injective arrows  of form $\{\bullet,\bullet'\}\lra X$.

We mentioned that the lifting property can be seen as a kind of negation. Confusingly, there are {\em two} negations, depending on whether the morphism appears on the left or right
side of the square, that are quite different: for example, both the pullback topology and the separation axiom $T_1$ are
negations of the same morphism, and the same goes for injectivity and injectivity on $\pi_0$ (see Figure~4(c,d)).

Now consider the standard example of something non-compact: the open covering
\[
\Bbb R= \bigcup\limits_{n\in\Bbb N} \{\,x\,:\, -n<x<n \,\}
\]
of the real line by infinitely many increasing intervals.
A related arrow in the category of topological spaces is
\[
\bigsqcup\limits_{n\in\Bbb N} \{\,x\,:\, -n<x<n \}\,\,\lra \Bbb R.
\]
Does the lifting property relative to that arrow define compactness? Not quite, but almost:
$$\{\}\lra X \rtt \bigsqcup\limits_{n\in\Bbb N} \{x\,:\, -n<x<n \,\}\,\,\lra \Bbb R$$
reads,  for $X$ connected,  as
``Every continuous real-valued function on $X$ is bounded, i.e. for each 
continuous $f:X\lra \Bbb R$ there is a natural number $n\in\Bbb N$
such that $-n<f(x)<n$ for each $x\in X$'', which is an early characterisation of compactness taught 
in a first course on analysis. Notice that this characterisation
mentions explicitly the arrow $X\lra\Bbb R$ and the bounded intervals of the real line, 
i.e. arrows $\{\,x\,:\, -n<x<n \}\xra{\,\subseteq\,} \Bbb R,\,{n\in\Bbb N}$
constituting the arrow-counterexample on the right hand side.

In a category of metric spaces with say distance non-increasing maps,
a metric space $X$ is {\em complete}, i.e. each Cauchy sequence $x_n\in X$, $n\in\Bbb N$, say
$dist(x_n,x_m)\leq 1/n$, converges to some point $x_\infty\in X$ such that $dist(x_\infty,x_n)\leq1/n$, iff 
\[
\left\{``x_n\text{''}\,:\ n\in{\Bbb N}\right\}\lra\left\{``x_n\text{''}\,:\,n\in {\Bbb N}\right\}\cup\{``x_\infty\text{''}\} \rtt X\lra\{\bullet\}
\]
(where $dist(``x_n\text{''},``x_m\text{''})=\frac1n$ for $m>n$, $dist(``x_\infty\text{''},``x_n\text{''})=\frac1n$, as defined above.)

In functional analysis, a (partially defined!) 
linear operator $f:X\lra Y$ between Banach spaces $X$ and $Y$ 
is {\em closed} iff for every convergent sequence $x_n\in X$, if $f(x_n)\xra[n\lra\infty]{}y$ in $Y$, 
then there is a $x\in X$ such that $f(x)=y$ and $x_n\xra[n\lra\infty]{}x$, i.e. 
\[
\left\{``x_n\text{''}\,:\ n\in{\Bbb N}\right\}\lra\left\{``x_n\text{''}\,:\,n\in {\Bbb N}\right\}\cup\{``x_\infty\text{''}\} \rtt \text{Domain}(f)\lra Y
\]

A module $P$ over a commutative ring $R$ is {\em projective} iff for an arbitrary arrow $N\lra M$ in the category of $R$-modules it holds
$$ 0\lra R \rtt N\lra M \ \implies\ 0\lra P \rtt N\lra M
.$$
Dually, a module $I$ over a ring $R$ is {\em injective} iff for an arbitrary arrow $N\lra M$ in the category of $R$-modules it holds
$$ R\lra 0 \rtt N\lra M \ \implies\  N\lra M \rtt I\lra 0\\
.$$
\def\nto{\not\to}
\def\ZpZ{{\Bbb Z}/\!p{\Bbb Z}}
\subsubsection{Finite groups.} There are examples outside of topology; see Appendix~\ref{app:examples}. 
Let us give some examples in group theory.
 There is no non-trivial homomorphism from a group $F$ to $G$, write $F\nto G$, iff  
$$0\lra F\,\rtt 0\lra G\text{ or equivalently }F\lra 0 \rtt G\lra 0.$$
A group $A$ is {\em Abelian} iff 
\[ \left<a,b\right> \,\lra\,\left<a,b:ab=ba\right>  \rtt\,\, A\lra 0
\]
where $\left<a,b\right> \,\lra\,\left<a,b:ab=ba\right>$ is the  abelianisation morphism sending the free group into the Abelian free group on two generators; 
a group $G$ is {\em perfect}, $G=[G,G]$, iff $G\nto A$ for any Abelian group $A$, i.e. 
\[ \left<a,b\right> \,\lra\,\left<a,b:ab=ba\right>  \rtt\,\, A\lra 0\ \implies\ G\lra 0 \rtt A\lra 0\]
in the category of finite or algebraic groups,
a group $H$ is {\em soluble} iff  $G\nto H$ for each perfect group $G$, 
i.e. 
$$0\lra G\,\rtt 0\lra H\text{ or equivalently }C\lra 0 \rtt H\lra 0.$$
A prime number $p$ does not divide the number  elements of a finite group $G$ 
iff  $G$ has no element of order $p$, i.e. no element $x\in G$ such that $x^p=1_G$ yet $x^1\neq 1_G,...,x^{p-1}\neq 1_G$,
equivalently $\ZpZ\nto G$, i.e.
$$0\lra \ZpZ\,\rtt 0\lra G\text{ or equivalently }\ZpZ\lra 0 \rtt G\lra 0.$$
A finite group $G$ is a $p$-group, i.e. the number of its elements is a power of a prime number $p$, iff in the category of finite groups 
$$0\lra \ZpZ \rtt 0\lra H \implies 0\lra H\rtt 0\lra G.$$

\subsection{Appendix. Iterated orthogonals: definitions and intuition.}

For a property (class) \ensuremath{C} of arrows (morphisms) in a category, 
define {\em its left and right orthogonals}, which we also call
{\em left and right negation}:
$$ C^l := \{ \ensuremath{f} :\text{ for each }g \in C\ \ensuremath{f} \,\rightthreetimes\,  \ensuremath{g} \} $$
$$ C^r := \{ \ensuremath{g} :\text{ for each }f \in C\ \ensuremath{f} \,\rightthreetimes\,  \ensuremath{g} \} $$
$$ C^{lr}:=(C^l)^r, ... $$
here $f \,\rightthreetimes\,  g$ reads ``$f$ has the left lifting property wrt $g$ '',
`` $f$ is (left) orthogonal to $g$ '',
i.e.  for  $f:A\longrightarrow B$, $g:X\longrightarrow Y$,
$f \,\rightthreetimes\, g$ iff for each $i:A\longrightarrow X$, $j:B\longrightarrow Y$ such that $ig=fj$ (``the square commutes''),
there is $j':B\longrightarrow X$ such that $fj'=i$ and $j'g=j$ (``there is a diagonal
making the diagram commute''), cf.~Fig.~\ref{fig:2a}.

The following observation is enough to reconstruct all the examples of
iterated orthogonals in this paper, with a bit of search and computation.
\begin{quote}
 {\bf Observation.}\\ 
  A number of elementary properties can be obtained by repeatedly passing
  to the left or right orthogonal $C^l, C^r, C^{lr}, C^{ll}, C^{rl}, C^{rr},...$
  starting from a simple class of morphisms, often
  a single (counter)example to the property you define.
  
 The counterexample is often implicit in the text of the definition of the property.
\end{quote}

A useful intuition is to think that the property of left-lifting against a
class \ensuremath{C} is a kind of negation of the property of being in \ensuremath{C}, and that
right-lifting is another kind of negation.  Hence the classes obtained from C
by taking orthogonals an odd number of times, such as $C^l, C^r, C^{lrl}, C^{lll}$
etc., represent various kinds of negation of \ensuremath{C}, so $C^l, C^r, C^{lrl}, C^{lll}$ each
consists of morphisms which are far from having property $C$.

Taking the orthogonal of a class \ensuremath{C} is a simple way to define a class of
morphisms excluding non-isomorphisms from \ensuremath{C}, in a way which is useful in a
diagram chasing computation.

The class $C^l$ is always closed under retracts, pullbacks, (small) products
(whenever they exist in the category) and composition of morphisms, and
contains all isomorphisms of $C$. Meanwhile, $C^r$ is closed under retracts,
pushouts, (small) coproducts and transfinite composition (filtered colimits) of
morphisms (whenever they exist in the category), and also contains all
isomorphisms.  Under some assumptions on existence of limits and colimits and ignoring 
set-theoretic difficulties\footnote{For an example of a theorem along these lines see~\href{https://core.ac.uk/download/pdf/82479252.pdf}{[Bousfield, Constructions of factorization systems in categories, 5.1 Ex, 3.1 Thm]}. Note that he considers the {\em unique} lifting property, unlike us.}, 
each morphism $X\lra Y$ decomposes both as 
$X\xra{(C)^l} \bullet \xra{(C)^{lr}} Y$ and $X\xra{(C)^{rl}} \bullet \xra{(C)^{r}} Y$.  

For example, the notion of isomorphism can be obtained starting from the class
of all morphisms, or any single example of an isomorphism:
$$
(Isomorphisms) = (all\ morphisms)^l = (all\ morphisms)^r = (h)^{lr} = (h)^{rl}
$$
where \ensuremath{h } is an arbitrary isomorphism.

\def\rrt#1#2#3#4#5#6{\xymatrix{ {#1} \ar[r]|{} \ar@{->}[d]|{#2} & {#4} \ar[d]|{#5} \\ {#3}  \ar[r] \ar@{-->}[ur]|{}& {#6} }}
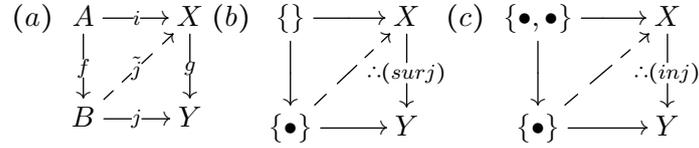
\begin{figure}
\begin{center}
$ (a)\ \xymatrix{ A \ar[r]|{i} \ar@{->}[d]|f & X \ar[d]|g \\ B \ar[r]|-{j} \ar@{-->}[ur]|{{\tilde j}}& Y }$
$(b)\  \rrt  {\{\}}  {} {\{\bullet\}}  X {\therefore(surj)} Y $
$(c)\  \rrt {\{\bullet,\bullet\}} {} {\{\bullet\}}  X {\therefore(inj)} Y $
\end{center}
\caption{\label{fig:2a}\normalsize
Lifting properties. 
 (a) The definition of a lifting property $f\rtt g$. 
 (b) $X\lra Y$ is surjective 
 (c) $X\lra Y$ is injective 
}
\end{figure}

\subsubsection{\label{app:rtt-examples}\label{app:examples}Examples of iterated orthogonals.}
Here give a list of examples of well-known properties which can be defined by 
iterated orthogonals starting from a simple class of morphisms.

\begin{itemize}
\item[ (i)] $(\emptyset\longrightarrow \{*\})^r$, $(0\longrightarrow R)^r$, and $\{0\longrightarrow \ZZ\}^r$ are the classes of surjections in
      in the categories of Sets, $R$-modules, and Groups, resp.,
      (where $\{*\}$ is the one-element set, and in the category of (not necessarily abelian) groups, $0$ denotes the trivial group)
\item[ (ii)] $(\{\star,\bullet\}\longrightarrow \{*\})^l=(\{\star,\bullet\}\longrightarrow \{*\})^r$, $(R\longrightarrow 0)^r$, $\{\ZZ\longrightarrow 0\}^r$ are  the classes of
      injections in the categories of Sets, $R$-modules, and Groups, resp
\item[ (iii)] in the category of $R$-modules, \\
                 \ \ a module \ensuremath{P} is projective iff $0\longrightarrow P$ is in $(0\longrightarrow R)^{rl}$\\
                 \ \  a module $I$ is injective iff $I\longrightarrow 0$ is in $(R\longrightarrow 0)^{rr}$
\item[ (iv)] in the category of Groups, \begin{itemize}
\item[]        a finite group \ensuremath{H} is nilpotent iff $H\longrightarrow H\times H$ is in $\{\, 0\longrightarrow \ensuremath{G} : G\text{ arbitrary} \}^{lr}$
\item[]        a finite group \ensuremath{H} is solvable iff $0\longrightarrow H$ is in $\{\, 0\longrightarrow A : A\text{ abelian }\}^{lr}= \{\, [G,G]\longrightarrow \ensuremath{G} : G\text{ arbitrary }\}^{lr}$
\item[]        a finite group \ensuremath{H} is a $p$-group iff $H\longrightarrow 0$ is in $\{\ZZ/p\ZZ\longrightarrow 0\}^{rr}$
\item[]        a group \ensuremath{F} is free iff $0\longrightarrow F$ is in  $\{0\longrightarrow \ZZ\}^{rl}$
\end{itemize}
\item[ (v)] in the category of metric spaces and uniformly continuous maps,\\
        a metric space \ensuremath{X} is complete iff $\{1/n\}_n\longrightarrow \{1/n\}_n\cup \{0\} \,\rightthreetimes\,  X\longrightarrow \{0\}$
           where the metric on $\{1/n\}_n$ and $\{1/n\}_n\cup \{0\}$ is induced from the
                 real line\\
        a subset $A \subset  X$ is closed iff  $\{1/n\}_n\longrightarrow \{1/n\}_n\cup \{0\} \,\rightthreetimes\,  A\longrightarrow X$
\item[ (vi)]    in the category of topological spaces,\\
     for a connected topological space $X$, each function on \ensuremath{X} is bounded
     iff 
         $$ \emptyset\longrightarrow \ensuremath{X} \,\rightthreetimes\,  \cup_n (-n,n) \longrightarrow  \RR$$

\item[ (vii)] in the category of topological spaces (see notation defined below),
\begin{itemize}
\item[] a space $X$ is path-connected iff $\{0,1\} \longrightarrow  [0,1] \,\rightthreetimes\,  \ensuremath{X} \longrightarrow  \{*\}$
\item[] a space $X$ is path-connected iff for each Hausdorff compact space $K$ and each injective map $\{x,y\} \hookrightarrow  K$ it holds
   $\{x,y\} \hookrightarrow  \ensuremath{K} \,\rightthreetimes\,  \ensuremath{X} \longrightarrow  \{*\}$
\end{itemize}
\end{itemize}
{\bf
Proof.
} In (iv), we use that a finite group $H$ is nilpotent iff the diagonal $\{
(h,h) : \ensuremath{h } \in \ensuremath{H} \}$ is subnormal in $H\times H$, 
cf.~\href{http://groupprops.subwiki.org/wiki/Nilpotent group}{[Nilp]}. 
\qed 

\subsection{A concise notation for certain properties in elementary point-set topology}

We introduce a concise, and in a sense intuitive, notation (syntax) able to express 
a number of properties in elementary point-set topology. It is appropriate 
for properties defined as iterated orthogonals (negation) starting from maps of finite topological spaces.\footnote
{I thank Urs Schreiber for help with the exposition in this subsection.}

For example, surjective, injective, connected, totally disconnected, and dense 
are expressed as 
$\{\{\}\lra \{a\}\}^r$, $\{\{x,y\}\lra \{x=y\}\}^r$, 
$\{\{x,y\}\lra \{x=y\}\}^l$ or $\{\{\}\lra \{a\}\}^{rll}$,
 $\{\{\}\lra \{a\}\}^{rllr}$,
$\{\{x\}\lra \{x\swarrow y\}\}^r$.

\subsubsection{\label{app:top-notation}Notation for maps between finite topological spaces.
}
A {\em topological space} comes with a {\em specialisation preorder} on its points: for
points $x,y \in X$,  $x \leq y$ iff $y \in cl x$ ($y$ is in the {\em topological closure} of $x$). 
The resulting {\em preordered set} may be regarded as a {\em category} whose
{\em objects} are the points of ${X}$ and where there is a unique {\em morphism} $x{\searrow}y$ iff $y \in cl x$.

For a {\em finite topological space} $X$, the specialisation preorder or
equivalently the corresponding category uniquely determines the space: a {\em
subset} of ${X}$ is {\em closed} iff it is
{\em downward closed}, or equivalently, 
is a full subcategory such that there are no morphisms going outside the subcategory.

The monotone maps (i.e. {\em functors}) are the {\em continuous maps} for this topology.

We denote a finite topological space by a list of the arrows (morphisms) in
the corresponding category; '$\leftrightarrow $' denotes an {\em isomorphism} and '$=$' denotes the {\em identity morphism}.  An arrow between two such lists denotes a {\em continuous map} (a functor) which sends each point to the correspondingly labelled point, but possibly turning some morphisms into identity morphisms, thus gluing some points. 

With this notation, we may display continuous functions for instance between the {\em discrete space} on two points, the {\em Sierpinski space}, the {\em antidiscrete space} and the {\em point space} as follows (where each point is understood to be mapped to the point of the same name in the next space):
$$
  \begin{array}{ccccccc}
  \{a,b\}
     &\longrightarrow& 
  \{a{\searrow}b\}
     &\longrightarrow& 
  \{a\leftrightarrow b\}
    &\longrightarrow& 
  \{a=b\}
  \\
  \text{(discrete space)}
     &\longrightarrow& 
  \text{(Sierpinski space)}
    &\longrightarrow& 
  \text{(antidiscrete space)}
    &\longrightarrow& 
  \text{(single point)}
  \end{array}
$$

In $A \longrightarrow  B$, each object and each morphism in $A$ necessarily appears in ${B}$ as well; we avoid listing 
the same object or morphism twice. Thus 
both 
$$
\{a\} \longrightarrow  \{a,b\}
  \phantom{AAA} \text{ and } \phantom{AAA} 
\{a\} \longrightarrow  \{b\}
$$ 
denote the same map from a single point to the discrete space with two points.
Both
 $$\{a{\swarrow}U{\searrow}x{\swarrow}V{\searrow}b\}\longrightarrow \{a{\swarrow}U=x=V{\searrow}b\}\text{ and }\{a{\swarrow}U{\searrow}x{\swarrow}V{\searrow}b\}\longrightarrow \{U=x=V\}$$
denote the morphism gluing points $U,x,V$.

In $\{a{\searrow}b\}$, the point $a$ is open and point ${b}$ is closed. We
denote points by $a,b,c,..,U,V,...,0,1..$ to make notation reflect the intended meaning, 
e.g.~$X\lra \{U\searrow U'\}$ reminds us that the preimage of $U$ determines an open subset of $X$, 
$\{x,y\}\lra X$ reminds us that the map determines points $x,y\in X$, and $\{o\searrow c\}$ reminds that $o$ is open and $c$ is closed. 

Each continuous map $A\lra B$ between finite spaces may be represented in this way; in the first list 
list relations between elements of $A$, and in the second list put relations between their images. 
However, note that this notation does not allow to represent {\em endomorphisms $A\lra A$}.
We think of this limitation 
as a feature and not a bug: in a diagram chasing computation, 
endomorphisms under transitive closure lead to infinite cycles, 
and thus our notation has better chance to define a computable fragment of topology.

\subsubsection{\label{app:rtt-top}Examples of iterated orthogonals obtained from maps between finite topological spaces.}
Here give a list of examples of well-known properties which can be defined by 
iterated orthogonals starting from maps between finite topological spaces, often with less than 5 elements.

In the category of topological spaces (see notation defined below),
\begin{itemize}
\item  a Hausdorff space $K$ is compact iff $K\longrightarrow \{*\}$ is in  $((\{o\}\longrightarrow \{o{\small\searrow}c\})^r_{<5})^{lr}$
\item  a  Hausdorff space $K$ is compact iff $K\longrightarrow \{*\}$ is in  $$
     \{\, \{a\leftrightarrow b\}\longrightarrow \{a=b\},\, \{o{\small\searrow}c\}\longrightarrow \{o=c\},\,
     \{c\}\longrightarrow \{o{\small\searrow}c\},\,\{a{\small\swarrow}o{\small\searrow}b\}\longrightarrow \{a=o=b\}\,\,\}^{lr}$$

\item  a space $D$ is discrete iff $ \emptyset \longrightarrow  D$ is in $   (\emptyset\longrightarrow \{*\})^{rl}      $

\item  a space $D$ is antidiscrete iff $ \ensuremath{D} \longrightarrow  \{*\} $ is in 
$(\{a,b\}\longrightarrow \{a=b\})^{rr}= (\{a\leftrightarrow b\}\longrightarrow \{a=b\})^{lr} $ 

\item  a space $K$ is connected or empty iff $K\longrightarrow \{*\}$ is in  $(\{a,b\}\longrightarrow \{a=b\})^l $
\item  a space $K$ is totally disconnected and non-empty iff $K\longrightarrow \{*\}$ is in  $(\{a,b\}\longrightarrow \{a=b\})^{lr} $

\item  a space $K$ is connected and non-empty
 iff
 for some arrow $\{*\}\longrightarrow K$\\
$\text{ \ \ \ \ \     } \{*\}\longrightarrow K$ is in
            $   (\emptyset\longrightarrow \{*\})^{rll} = (\{a\}\longrightarrow \{a,b\})^l$

\item  a space $K$ is non-empty iff $K\longrightarrow \{*\}$ is in $   (\emptyset\longrightarrow \{*\})^l$
\item  a space $K$ is empty iff $K \longrightarrow \{*\}$ is in $   (\emptyset\longrightarrow \{*\})^{ll}$
\item  a space $K$ is $T_0$ iff $K  \longrightarrow \{*\}$ is in $   (\{a\leftrightarrow b\}\longrightarrow \{a=b\})^r$
\item   a space $K$ is $T_1$ iff $K  \longrightarrow \{*\}$ is in $   (\{a{\small\searrow}b\}\longrightarrow \{a=b\})^r$
\item  a space $X$ is Hausdorff iff for each injective map $\{x,y\} \hookrightarrow  X$
it holds $\{x,y\} \hookrightarrow  \ensuremath{X} \,\rightthreetimes\,  \{ \ensuremath{x} {\small\searrow} \ensuremath{o} {\small\swarrow} \ensuremath{y} \} \longrightarrow  \{ x=o=y \}$

\item  a non-empty space $X$ is regular (T3) iff for each arrow $    \{x\} \longrightarrow  X$ it holds
    $    \{x\} \longrightarrow  \ensuremath{X} \,\rightthreetimes\,  \{x{\small\searrow}X{\small\swarrow}U{\small\searrow}F\} \longrightarrow  \{x=X=U{\small\searrow}F\}$
\item  a space $X$ is normal (T4) iff $\emptyset \longrightarrow \ensuremath{X} \,\rightthreetimes\,   \{a{\small\swarrow}U{\small\searrow}x{\small\swarrow}V{\small\searrow}b\}\longrightarrow \{a{\small\swarrow}U=x=V{\small\searrow}b\}$

\item  a space $X$ is completely normal iff $\emptyset\longrightarrow \ensuremath{X} \,\rightthreetimes\,  [0,1]\longrightarrow \{0{\small\swarrow}x{\small\searrow}1\}$
 where the map $[0,1]\longrightarrow \{0{\small\swarrow}x{\small\searrow}1\}$ sends $0$ to $0$, $1$ to $1$, and the rest $(0,1)$ to $x$
\item  a space $X$ is path-connected iff $\{0,1\} \longrightarrow  [0,1] \,\rightthreetimes\,  \ensuremath{X} \longrightarrow  \{*\}$
\item  a space $X$ is path-connected iff for each Hausdorff compact space $K$ and each injective map $\{x,y\} \hookrightarrow  K$ it holds
   $\{x,y\} \hookrightarrow  \ensuremath{K} \,\rightthreetimes\,  \ensuremath{X} \longrightarrow  \{*\}$

\item  a non-empty space $X$ is regular (T3) iff for each arrow $    \{x\} \longrightarrow  X$ it holds
    $    \{x\} \longrightarrow  \ensuremath{X} \,\rightthreetimes\,  \{x{\small\searrow}X{\small\swarrow}U{\small\searrow}F\} \longrightarrow  \{x=X=U{\small\searrow}F\}$
\item  a space $X$ is normal (T4) iff $\emptyset \longrightarrow \ensuremath{X} \,\rightthreetimes\,   \{a{\small\swarrow}U{\small\searrow}x{\small\swarrow}V{\small\searrow}b\}\longrightarrow \{a{\small\swarrow}U=x=V{\small\searrow}b\}$

\item  a space $X$ is completely normal iff $\emptyset\longrightarrow \ensuremath{X} \,\rightthreetimes\,  [0,1]\longrightarrow \{0{\small\swarrow}x{\small\searrow}1\}$
 where the map $[0,1]\longrightarrow \{0{\small\swarrow}x{\small\searrow}1\}$ sends $0$ to $0$, $1$ to $1$, and the rest $(0,1)$ to $x$
\item  a space $X$ is path-connected iff $\{0,1\} \longrightarrow  [0,1] \,\rightthreetimes\,  \ensuremath{X} \longrightarrow  \{*\}$
\item  a space $X$ is path-connected iff for each Hausdorff compact space $K$ and each injective map $\{x,y\} \hookrightarrow  K$ it holds
   $\{x,y\} \hookrightarrow  \ensuremath{K} \,\rightthreetimes\,  \ensuremath{X} \longrightarrow  \{*\}$

\item        $(\emptyset\longrightarrow \{*\})^r$   is the class of surjections
\item        $(\emptyset\longrightarrow \{*\})^{rr}$ is the class of subsets, i.e. injective maps $A\hookrightarrow B$ where the topology on $A$ is induced from $B$
\item        $(\emptyset\longrightarrow \{*\})^{lll}$ is the class of maps $A\longrightarrow B$ which split

\item        $(\{b\}\longrightarrow \{a{\small\searrow}b\})^l$ is the class of maps with dense image
\item        $(\{b\}\longrightarrow \{a{\small\searrow}b\})^{lr}$ is the class of closed subsets $A \subset  X$, $A$ a closed subset of $X$
\item        $((\{a\}\longrightarrow \{a{\small\searrow}b\})^r_{<5})^{lr}$ is roughly the class of proper maps
       (see below).
\end{itemize}
{\bf
Proof.
} 
Items related to compactness and proper maps are discussed in ??. 
Other items require a simple if tedious verification.  \qed

\subsection{Separation axioms as orthogonals.} 
See \url{https://ncatlab.org/nlab/show/separation+axioms+in+terms+of+lifting+properties} for a list of reformulations of the separation axioms.

\subsection{\label{app:AEEA}Appendix. Compactness as being uniform: change of order of quantifiers}

We give several examples where an application of compactness
can be reformulated as changing the order of quantifiers in a formula.

\subsubsection{ 
Each real-valued function on a compact set is bounded
}
$$\frac{
\forall x \in K \exists M  ( f(x) < M )
}{
\exists M \forall x \in K  ( f(x) < M )
}$$
%
Note this is a lifting property, for $K$ connected:

$$\{\}\longrightarrow  K \rtt   \sqcup_{ n\in\NN} (-n,n) \longrightarrow   \RR$$ 
here  $\sqcup_ n (-n,n) \longrightarrow   \RR$ denotes the map to the real line
from the disjoint union of intervals $(-n,n)$ which cover it. 
Note this is a standard example of an open covering of $\RR$ which 
shows it is not compact.

\subsubsection{ 
The image of a closed set is closed
} 

$K$ is compact iff the following implication  
holds for each set $X$ and each subset $Z \subset \XxK$:
$$\frac{
\forall y \in K \exists U \exists V ( U \subset X\text{ open and }V \subset K\text{ open and }a \in U\text{ and }y \in V\text{ and }\UxV \subset Z )
}{
\exists U \exists V \forall y \in K  ( U \subset X\text{ open and }V \subset K\text{ open and }a \in U\text{ and }y \in V\text{ and }\UxV \subset Z )
}$$
%
%
%
%
%
The hypothesis says $Z$ contains a rectangular open neighbourhood of each point of the line $\{a\}\times K$;
the conclusion says that $Z$ contains a rectangular open neighbourhood of the whole line $\{a\}\times K$.

\subsubsection{ 
A Hausdorff compact is necessarily normal.
} 

The application of compactness in the usual proof of this implication 
amounts to the following change of order of quantifiers:
$$\frac{
\forall a \in A \forall b \in B \exists U \exists V ( a \in U\text{  and }b \in V\text{ and }U\cap V=\{\}\text{ and }U \subset K\text{ open and }V \subset K\text{ open} )
}{
\exists U \exists V \forall a \in A \forall b \in B  ( a \in U\text{  and }b \in V\text{ and }U\cap V=\{\}\text{ and }U \subset K\text{ open and }V \subset K\text{ open} )}$$
%
%
\subsubsection{ 
Lebesgue's number Lemma 
} 

Let $S$ be a family of (arbitrary) subsets of a metric space $X$.
$$\frac{
\forall x \in X \exists \delta>0 \exists U \in S \forall y \in X  ( dist(x,y)<\delta  \implies y \in U )
}{
\exists \delta>0 \forall x \in X \exists U \in S \forall y \in X  ( dist(x,y)<\delta \implies y \in U )
}$$
%
%
The hypothesis says that $\{ Inn\, U\ :\ U \in S \}$ is an open cover of $X$; 
the conclusion is as usually stated, that each set of diameter $<\delta$ is covered by a single member of the cover. 

Note that this lemma may be expressed in terms of uniform structures.

\subsubsection{Paracompactness.}

Alexandroff writes ``as it seems to me, one of the deepest and most interesting 
properties of paracompacts'' is the following theorem of A.Stone: 
that \bqqq
  A $T_1$-space is {\em paracompact} iff for each open covering $\alpha$ of $X$ there is an open covering $\beta$
such that for each $x$ in $X$ there is $U$ in $A$ such that
         $\cup \{ V \in B : x \in V \} \subset U $   
\eqqq
As quantifier exchange, this is: 
$$\frac{
\text{for each open covering }\alpha\text{ exists open covering }\beta.\ \forall x\in X\forall V \in \beta \exists U \in \alpha ( x \in V \implies V \subset U )
}{
\text{for each open covering }\alpha\text{ exists open covering }\beta.\ \forall x\in X\exists U\in\alpha \forall V \in \beta ( x \in V \implies V \subset U )
}$$
%
%
%
The hypothesis holds trivially: take $\beta=\alpha , V=U$. 

\begin{enonce}{Question} 
Describe a logic and a class of formulae where such exchange of order 
quantifiers is permissible. Is there a treatment of compactness in terms
of changing order of quantifiers ? 
\end{enonce} 

\subsection*{Acknowledgements.} To be written.

 This work is a continuation of [DMG]; early history is given there. I thank
M.Bays, G.Cherlin, D.Krachun, K.Pimenov, V.Sosnilo, S.Sinchuk and P.Zusmanovich for discussions and proofreading; I
thank L.Beklemishev, N.Durov, S.V.Ivanov, D.Grayson, S.Podkorytov, A.L.Smirnov for discussions. I also
thank several students for encouraging and helpful discussions.  Chebyshev
laboratory, St.Petersburg State University, provided a coffee machine and an
excellent company around it to chat about mathematics.  Special thanks are to
Martin Bays for many corrections and helpful discussions. Several observations
in this paper are due to Martin Bays. I thank S.V.Ivanov for several
encouraging and useful discussions; in particular, he suggested to look at the
Lebesque's number lemma and the Arzela-Ascoli theorem. A discussion with Sergei
Kryzhevich motivated the group theory examples.

   Much of this paper was done in St.Petersburg; it wouldn't have been possible
without support of family and friends who created an excellent social
environment and who occasionally accepted an invitation for a walk or a coffee
or extended an invitation; alas, I made such a poor use of it all.

  This note is elementary, and it was embarrassing and boring, and
embarrassingly boring, to think or talk about matters so trivial, but luckily
I had no obligations for a time.

\end{document}